\documentclass[final,hidelinks,onefignum,onetabnum]{siamart220329}

\usepackage{lipsum}
\usepackage{amsfonts}
\usepackage{graphicx}
\usepackage{epstopdf}
\usepackage{algorithmic}
\usepackage{algorithm}

 \usepackage{booktabs, ragged2e, array, adjustbox}
 \usepackage{array} 

\usepackage{chngcntr}

\counterwithout{algorithm}{section}

\usepackage{CJKutf8}

\usepackage{relsize} 

\usepackage{graphicx} 

\usepackage{enumitem}
\usepackage{amsmath} 
\usepackage{amssymb}
\usepackage{mathtools}
\usepackage{esint}
\usepackage{mathrsfs}
\ifpdf
  \DeclareGraphicsExtensions{.eps,.pdf,.png,.jpg}
\else
  \DeclareGraphicsExtensions{.eps}
\fi

\usepackage{tabularx}

\usepackage{graphicx}
\usepackage{subcaption}

\newsiamremark{defintion}{Defintion}
\newsiamremark{remark}{Remark}
\newsiamremark{hypothesis}{Hypothesis}
\newtheorem{example}{Example}
\crefname{hypothesis}{Hypothesis}{Hypotheses}
\newsiamthm{claim}{Claim}

\title{An Iterative Method with Asymptotic Orthogonality for Simultaneous Eigenpair Computation\thanks{Submitted to the editors DATE.
		\funding{This work was supported by the National Natural Science Foundation of China under grant 12571446 and the National Key R \& D Program of China under grants 2025YFA1016600 and 2025YFA1016601.}}}

\author{Shengyue Wang\thanks{State Key Laboratory of Mathematical Sciences (SKLMS), Academy of Mathematics and Systems Science, Chinese Academy of Sciences, Beijing 100190, China; and School of Mathematical Sciences, University of Chinese Academy of Sciences, Beijing 100049, China
 (\email{wangshengyue@amss.ac.cn}, \email{azhou@lsec.cc.ac.cn}).
}\and Aihui Zhou\footnotemark[2]}

\usepackage{amsopn}

\usepackage{cleveref}

\ifpdf
\hypersetup{
	pdftitle={An Iterative Method with Asymptotic Orthogonality for Simultaneous Eigenpair Computation},
	pdfauthor={S.Wang and A.Zhou}
}
\fi

\begin{document}
	
	\maketitle
	
	\tableofcontents 

	\begin{abstract}
		The simultaneous computation of a cluster of eigenpairs with mutually orthogonal eigenvectors is a basic task in scientific computing. We develop a predictor--corrector discretization of the quasi-Grassmannian gradient flow for simultaneous eigenpair computation. The proposed iteration requires neither orthogonal initial data nor any orthogonalization operation: the predictor preserves the current Gram matrix, whereas the corrector reduces the orthogonality error, so that the iterates approach orthogonality asymptotically. We establish well-posedness of the discrete scheme and prove invariant-subspace nonexpansion, asymptotic orthogonality, energy decrease, and convergence to the target eigenspace. Numerical experiments with the discrete equations solved approximately show that the numerical iterates converge to eigenpair approximations whose eigenvectors are numerically orthogonal to high accuracy, while the energy and gradient norm decrease over the iterations.
	\end{abstract}
	
	

	\section{Introduction}
	
	Computing many eigenpairs simultaneously arises in scientific computing. Many applications require not only the eigenvalues but also a mutually orthonormal family of corresponding eigenvectors. In electronic structure calculations, the occupied states satisfy orthonormality constraints \cite{le2005computational,leszczynski2012handbook, martin2020electronic}. In structural vibration and photonic crystal computations, multiple eigenmodes are required to construct modal and spectral representations \cite{bathe1973solution,dorfler2011photonic}. In addition, an orthonormal basis of the target invariant subspace enables stable spectral projection and prevents redundant representations of the same subspace \cite{golub2013matrix,saad1992numerical}. These considerations make the simultaneous computation of many orthogonal eigenvectors relevant to large-scale applications. This paper focuses on eigenvalue problems associated with self-adjoint partial differential operators. After discretization, the method and analysis below are formulated in the resulting finite-dimensional space.

	For such block eigenvalue computations, maintaining mutual orthogonality may represent part of the computational cost. Block Krylov subspace methods, Davidson-type methods, and subspace iterations commonly employ periodic orthogonalization and Rayleigh--Ritz procedures \cite{davidson1975iterative,golub2013matrix,saad1992numerical}. Here, orthogonalization means replacing an iterate by an orthonormal basis through a QR factorization, a Gram--Schmidt procedure, or an equivalent retraction. These operations require additional work when many eigenvectors are sought and may involve repeated synchronization in distributed-memory computations \cite{demmel2012communication,leon2013gram}. Polynomially filtered subspace iterations, in particular Chebyshev-filtered methods, can perform several filtering steps between orthogonalizations and thereby reduce their frequency and relative cost \cite{winkelmann2019chase,zhou2006chebyshev}. Nevertheless, avoiding repeated basis orthogonalization remains relevant in large-scale block computations. This motivates the design of iterative schemes that produce mutually orthogonal eigenvector approximations without any additional orthogonalization.

	One way to construct such iterations is to characterize the target eigenspace as an energy minimizer under orthogonality constraints. This leads to gradient flows on the Stiefel or Grassmann manifold. Dai et al. introduced a gradient-flow model for the Kohn--Sham energy and a midpoint discretization that preserves orthogonality and converges to a minimizer in a neighborhood of the solution \cite{dai2020}. They subsequently developed a class of orthogonality-preserving time discretizations, including midpoint- and Crank--Nicolson-type schemes, together with convergence analysis and an adaptive time-step strategy \cite{dai2021convergent}. For linear eigenvalue problems, Chu et al. proposed an orthogonality-preserving evolution model and an explicit time discretization that dissipates energy and can rapidly produce eigenpair approximations \cite{chu2025orthogonality}. When initialized on the Stiefel manifold, these schemes preserve orthogonality in arithmetic without roundoff. In finite-precision computation, orthogonality is generally attained only up to a prescribed numerical tolerance. {The skew-symmetric evolutions underlying these methods preserve the Gram matrix and hence leave the orthogonality error of nonorthogonal data unchanged. To accommodate nonorthogonal initial data in the simultaneous computation of many eigenpairs, a quasi-Grassmannian gradient flow was introduced in \cite{Wang2025}. This flow drives the Gram matrix to the identity and thereby yields asymptotic orthogonality without orthogonalization.}

	{Building on this continuous model, we propose a predictor--corrector iteration with asymptotic orthogonality.} Given the current iterate \(U_n\), the predictor is a midpoint discretization generated by a skew-symmetric operator and preserves its Gram matrix, while the corrector acts directly on the orthogonality error matrix \(I_N-\langle U_n,U_n\rangle\). For the nonorthogonal initial data considered below, the proposed scheme preserves full column rank and satisfies
	\begin{equation*}
		0<\langle U_n,U_n\rangle\leqslant I_N,
		\qquad
		\lim_{n\rightarrow\infty}
		\left\|I_N-\langle U_n,U_n\rangle\right\|=0.
	\end{equation*}
	For orthonormal initial data, the scheme preserves orthogonality. It does not require any additional QR, Gram--Schmidt, or manifold-retraction step.

	The main contributions of this work are summarized as follows:
	\begin{itemize}
		\item {We construct a predictor--corrector discretization of the quasi-Grassmannian gradient flow \cite{Wang2025} for the simultaneous approximation of a cluster of eigenpairs.} {The predictor, formed from the gradient of the associated energy, preserves the current Gram matrix, whereas the corrector moves it toward the identity and thereby reduces the orthogonality error.} The scheme permits nonorthogonal initial data and requires no additional orthogonalization.

		\item We establish {existence and uniqueness of solutions to \eqref{Discretization scheme}} near the target eigenspace. For the resulting iteration, we prove {invariant-subspace nonexpansion and preservation of full column rank}, together with asymptotic orthogonality and an exponential decay estimate for the orthogonality error. We further establish energy decay and convergence to the target eigenspace from a neighborhood of the solution.

		\item We report numerical experiments for three eigenvalue problems. In these tests, the energy, gradient norm, and difference from the final computed iterate decrease, while the orthogonality error reaches the prescribed tolerance without any orthogonalization.
	\end{itemize}

	The remainder of the paper is organized as follows. Section \ref{sec:Problem setting} introduces the continuous and discretized eigenvalue problems together with the notation used in the analysis. Section \ref{sec:asymptotic-orthogonality algorithm} presents the predictor--corrector scheme and establishes {invariant-subspace nonexpansion} and asymptotic orthogonality. Section \ref{sec:Numerical analysis} proves the solvability, energy decay, and convergence results. Section \ref{Numerical experiments} presents the practical iteration and reports the numerical results, which show that the proposed scheme has the potential to compute a cluster of eigenpairs simultaneously. Section \ref{sec:Concluding remarks} concludes the paper. The appendices provide the proofs of four auxiliary convergence results and the derivation and implementation details of the practical iteration.

	\section{Problem setting} \label{sec:Problem setting}
	{We first formulate the continuous and discretized eigenvalue problems, relate them to the quasi-Grassmannian gradient flow, and then introduce the notation used in the analysis.}

	Let $\Omega \subset \mathbb{R}^d$ $(d\in \mathbb{N}_+)$ be a bounded domain with a regular boundary (see, e.g. \cite{evans2022partial,zeidler2013nonlinear}). We denote $H^k(\Omega)$, $k \geqslant 0$ the standard Sobolev space, and set
	$$
	L^2(\Omega)=H^0(\Omega) \quad \text { and } \quad H_0^1(\Omega)=\left\{u\in H^1(\Omega):u|_{\partial \Omega}=0  \right\},
	$$
	where $u|_{\partial \Omega}=0 $ is understood in the sense of trace. We define the inner product $(\cdot, \cdot)$ and norm $\|\cdot\|$ of the space $L^2(\Omega)$ respectively as
	$$
	(u, v)=\int_{\Omega} u v \quad \text { and } \quad\|u\|=\sqrt{(u, u)}.
	$$
	{Let $V_{\text{ext}}:\Omega\rightarrow\mathbb R$ be a potential such that the symmetric bilinear form
	\begin{equation*}
		a(u,v)=\frac{1}{2}(\nabla u,\nabla v)+(V_{\text{ext}}u,v),
		\qquad u,v\in H_0^1(\Omega),
	\end{equation*}
	is continuous and satisfies the G{\aa}rding inequality. That is, there exist $c_0,c_1>0$ and $c_2\geqslant0$ such that
	\begin{equation*}
		|a(u,v)|\leqslant c_0\|\nabla u\|\|\nabla v\|,
		\qquad
		a(u,u)\geqslant c_1\|\nabla u\|^2-c_2\|u\|^2
	\end{equation*}
	for all $u,v\in H_0^1(\Omega)$. Let $\mathcal H=-\frac{1}{2}\Delta+V_{\text{ext}}$ denote the self-adjoint operator associated with $a$, characterized by
	\begin{equation*}
		\langle\mathcal Hu,v\rangle=a(u,v)
		\qquad \forall u,v\in H_0^1(\Omega),
	\end{equation*}
	where $\langle\cdot,\cdot\rangle$ denotes the dual pairing between $H^{-1}(\Omega)$ and $H_0^1(\Omega)$.} For the subsequent analysis, we further assume that \(\mathcal{H}\) has at least \(N\) negative eigenvalues.
	
	We seek the $N$ smallest eigenvalues $\lambda_1\leqslant \lambda_2\leqslant \cdots\leqslant \lambda_N    \ (<\lambda_{N+1})$ and corresponding eigenfunctions \(\{\tilde{v}_j\}_{j=1}^N\) of \(\mathcal{H}\), characterized by the eigenvalue problem:
	\begin{equation}
		\label{1-dim linear eigenvalue problem}
		\left\{\begin{aligned}
			&\langle \mathcal{H}  \tilde{v}_j , v\rangle =\lambda_j\langle {   \tilde{v}_j},v\rangle, \ \forall v\in H_0^1(\Omega)\quad  j=1,2,\cdots,N,
			\\&({   \tilde{v}_i} ,{   \tilde{v}_j} ) =\delta_{ij},
			\quad i,j=1,2,\ldots,N,
		\end{aligned}\right.
	\end{equation}
	where \(\delta_{ij}\) is the Kronecker delta.
	
	Consider the product Hilbert space
	\begin{equation*}
			\big(H_0^1(\Omega)\big)^N = \left\{(u_1,u_2, \cdots,u_N): u_i \in H_0^1(\Omega), \ i = 1,2,\cdots,N\right\},
	\end{equation*}
		equipped with the inner product matrix
		\begin{equation*}
				\langle U,V\rangle = \left(\left(u_i, v_j\right)\right)_{i,j = 1}^N \quad \forall U,V \in \big(H_0^1(\Omega)\big)^N,
		\end{equation*}
		and the associated inner product and norm
		\begin{equation*}
				\left(U,V\right) =\operatorname{tr}\left(\langle U,V\rangle\right),	\qquad \left\| U\right\| = \left(U,U\right)^{\frac{1}{2}} \quad \forall U,V \in \big(H_0^1(\Omega)\big)^N.  
		\end{equation*}
		For $v\in H_0^1(\Omega)$ and $U=(u_1,u_2,\ldots,u_N)\in\big(H_0^1(\Omega)\big)^N$, we also write
		\begin{equation*}
			\langle v,U\rangle=\left((v,u_j)\right)_{j=1}^N.
		\end{equation*}
		For $F=(f_1,f_2,\ldots,f_N)\in\big(H^{-1}(\Omega)\big)^N$, the matrix-valued dual pairings are defined by
		\begin{equation*}
			\langle F,U\rangle=\left(\langle f_i,u_j\rangle\right)_{i,j=1}^N.
		\end{equation*}
	
	We introduce the Stiefel manifold as follows
		$$	\mathcal{M}^N=\left\{U \in \big(H_0^1(\Omega)\big)^N: \langle U,U\rangle=I_N\right\},$$
	where \(I_N\) denotes the \(N \times N\) identity matrix, and \(\langle U,U\rangle=((u_i,u_j))_{i,j=1}^N\) is the Gram matrix of \(U\). Notably, it holds that
	\begin{equation*}
			U \in \mathcal{M}^N \iff UQ \in \mathcal{M}^N \quad \forall Q \in \mathcal{O}^N.
	\end{equation*}
	Here, \(\mathcal{O}^N\) denotes the set of \(N \times N\) orthogonal matrices.
	
	We extend the operator $\mathcal{H}$ to $\big(H_0^1(\Omega)\big)^N$, also denoted as $\mathcal{H}$ for simplicity. That is, for any $U=\left(u_1, u_2, \cdots, u_N\right) \in\left(H_0^1(\Omega)\right)^N$,
	\begin{equation*}
			\mathcal{H}U =\left(\mathcal{H}u_1,\mathcal{H} u_2, \cdots, \mathcal{H}u_N\right). 
	\end{equation*} 
	Under this definition, (\ref{1-dim linear eigenvalue problem}) is equivalent to
	\begin{equation}
		\label{whole space linear eigenvalue problem}
		\left\{ \begin{aligned}
			&\langle \mathcal{H} {   \tilde{V}_{*}} ,U\rangle = \langle {   \tilde{V}_{*}} \Lambda, U \rangle \quad \forall U \in \big(H_0^1(\Omega)\big)^N,     
			\\& {   \tilde{V}_{*}}   \in \mathcal{M}^N,
		\end{aligned}\right.
	\end{equation}
	where $\Lambda =\operatorname{diag}(\lambda_1,\lambda_2,\cdots, \lambda_N)$ and the columns of ${   \tilde{V}_*=\left(\tilde{v}_1, \tilde{v}_2, \ldots, \tilde{v}_N\right)} \in\left(H_0^1\left(\Omega\right)\right)^N$ are the corresponding eigenfunctions.

	Let $\mathcal{G}^N$ denote the Grassmann manifold (also known as Grassmannian), a quotient manifold of Stiefel manifold defined by
	$$\mathcal{G}^N = \mathcal{M}^N/ \sim , $$
	where $\sim$ is the equivalence relation defined as: $\hat{U} \sim U$ if and only if $\hat{U}=U Q$ for some $Q \in \mathcal{O}^N$. For any $U\in \big(H_0^1(\Omega)\big)^N$, the equivalence class is denoted by
	$$
	[U]=\left\{U Q: Q \in \mathcal{O}^N\right\},
	$$
 and \(\mathcal{G}^N\) is thus defined as
	$$
	\mathcal{G}^N=\left\{[U]: U \in \mathcal{M}^N\right\} .
	$$

	We define the energy by
	\begin{equation*}
		E(U)\stackrel{\Delta}{=}
		\frac{1}{2}\operatorname{tr}\left(\langle U,\mathcal{H}U\rangle\right),
		\qquad U\in\big(H_0^1(\Omega)\big)^N.
	\end{equation*}
	It is invariant under right multiplication by an orthogonal matrix:
	\begin{equation*}
		E(UQ)=E(U)
		\qquad \forall U\in\big(H_0^1(\Omega)\big)^N,
		\quad \forall Q\in\mathcal O^N.
	\end{equation*}
	Thus, \eqref{whole space linear eigenvalue problem} can be formulated as
	\begin{equation}
		\label{whole space minimization in Grassmann}
		\min_{[U]\in\mathcal G^N} E(U).
	\end{equation}
	Moreover, $[\tilde{V}_*]$ is the unique minimizer of (\ref{whole space minimization in Grassmann}) provided the spectral gap \(\lambda_{N+1} - \lambda_N > 0\) (see \cite{schneider2009direct}).

	In practice, \eqref{whole space linear eigenvalue problem} admits various discretization methods, including the plane wave method, the finite difference method, and the finite element method. In this paper, we work within the spatially discretized space $\left(\mathcal{V}^{N_g}\right)^N,$ where $\mathcal{V}^{N_g} \subset H_0^1(\Omega)$ is an $N_g$-dimensional subspace ($N_g \gg N$ typically). Within $\left(\mathcal{V}^{N_g}\right)^N$, we define the discretized Stiefel manifold as
	\begin{equation*}
			\mathcal{M}^{N;N_g}=\left\{U \in  \left(\mathcal{V}^{N_g}\right)^N : \langle U,U\rangle=I_N\right\},
	\end{equation*}
	with the corresponding discretized quotient manifold (the discretized Grassmann manifold):
	\begin{equation*}
			\mathcal{G}^{N;N_g} =	\mathcal{M}^{N;N_g} / \sim.
	\end{equation*}
	
	Denote \(H= \mathcal{P}_{\mathcal{V}^{N_g}}\mathcal{H}\mathcal{P}_{\mathcal{V}^{N_g}}:(\mathcal{V}^{N_g})^N\rightarrow(\mathcal{V}^{N_g})^N\), where \( \mathcal{P}_{\mathcal{V}^{N_g}}\) is the orthogonal projection onto \((\mathcal{V}^{N_g})^N\). We denote by \(V_*=\left(v_1, v_2, \cdots, v_N\right)\in\mathcal{M}^{N;N_g}\) the orthonormal eigenvectors associated with the first \(N\) discrete eigenvalues of \(H\). With a slight abuse of notation, these discrete eigenvalues are still denoted by \(\lambda_1,\lambda_2,\cdots,\lambda_N\). The associated discretized eigenvalue problem and minimization problem can be formulated as, respectively,
	\begin{equation}	\label{linear eigenvalue problem}
			\left\{ \begin{aligned}
			&H V_{*} = V_{*}\Lambda,     
			\\& V_{*}\in \mathcal{M}^{N;N_g} ,
		\end{aligned}\right.
	\end{equation}
	\begin{equation}
	\label{minimization in Grassmann}
	\min _{[U] \in \mathcal{G}^{N;N_g}} E (U),
	\end{equation}
	where $E(U) = \frac{1}{2}\operatorname{tr}\left(\langle U,\mathcal{H}U\rangle\right) =  \frac{1}{2}\operatorname{tr}\left(\langle U,HU\rangle\right), \forall U \in(\mathcal{V}^{N_g})^N $ and $[V_{*}]$ is the unique minimizer of (\ref{minimization in Grassmann}) under the discrete spectral gap \(\lambda_N<\lambda_{N+1}\).

	The gradient of \(E\) in \((\mathcal{V}^{N_g})^N\) is
	\begin{equation*}
		\nabla E(U)=HU.
	\end{equation*}
	For \(U\in\mathcal{M}^{N;N_g}\), the Grassmannian gradient of \(E\) at \([U]\) is
	\begin{equation*}
		\nabla_G E(U)
		=
		\nabla E(U)-U\langle U,\nabla E(U)\rangle
		=
		HU-U\langle U,HU\rangle.
	\end{equation*}
	{Following \cite{Wang2025}, for \(U\notin\mathcal{M}^{N;N_g}\), we use the same expression as an extension of the Grassmannian gradient to \((\mathcal{V}^{N_g})^N\).} For brevity, \(\nabla_G E(U)\) is referred to as the gradient below. In particular, \(\nabla_G E(V_*)=0\). {The quasi-Grassmannian gradient flow introduced in \cite{Wang2025} extends the Grassmannian gradient to the ambient space. In the infinite-dimensional operator setting, \cite{Wang2025} establishes the existence and uniqueness of a weak solution and derives its analytic representation. Its quasi-orthogonal character is reflected in the following properties:
	\begin{itemize}
		\item The initial data need not lie on the Stiefel manifold, and the solution is not required to satisfy the orthogonality constraint at each time.
		\item For the nonorthogonal initial data considered in \cite{Wang2025}, the Gram matrix converges to the identity. Hence, orthogonality is attained asymptotically through the evolution without orthogonalization.
		\item The gradient converges to zero, and the equivalence class of the solution converges to the target eigenspace.
	\end{itemize}
	After spatial discretization, the flow associated with \eqref{minimization in Grassmann} takes the form
	\begin{equation}
		\label{eq: continuous model}
		\left\{\begin{aligned}
			\frac{\mathrm{d}U}{\mathrm{d}t}&=-\nabla_G E(U), && t>0,\\
			U(0)&=U_0, && U_0\in(\mathcal{V}^{N_g})^N.
		\end{aligned}\right.
	\end{equation}
	Every representative \(V_*Q\), \(Q\in\mathcal O^N\), of the target eigenspace is a stationary solution of this flow, while \([V_*]\) is characterized equivalently by \eqref{linear eigenvalue problem} and \eqref{minimization in Grassmann}.}

	To facilitate subsequent discussions, we summarize key definitions used throughout the paper in Table \ref{table:notation}. All symbols are consistently applied in theoretical analyses, algorithm constructions, and numerical experiments.
	\begin{table}[htbp]
		\centering
		\caption{Notation}
		\label{table:notation}
		\adjustbox{max width=.85\linewidth}{%
			\begin{tabular}{>{\RaggedRight\arraybackslash}p{0.25\linewidth}
					>{\RaggedRight\arraybackslash}p{0.69\linewidth}}
				\toprule[1.2pt]
				\textbf{Symbol} & \textbf{Definition \& Explanation} \\
				\midrule
				
				$A \leqslant B$ &
				Positive semidefinite partial order for symmetric matrices
				$A,B\in\mathbb{R}^{N\times N}$:
				\par\smallskip
				\centerline{$
					A\leqslant B\iff a^{\!\top}\!Aa\leqslant a^{\!\top}\!Ba
					\;\;\forall a\in\mathbb{R}^{N}
					$}\\[2pt]
				
				$\lambda_{\max}(A)$ / $\lambda_{\min}(A)$ &
				Largest / smallest eigenvalue of matrix $A$\\[10pt]
				
				$\sigma_{\max}(U)$ &
				Maximum singular value of
				$U\in(\mathcal{V}^{N_{g}})^{N}$:
				\par\smallskip
				\centerline{$
					\sigma_{\max}(U)=\sqrt{\lambda_{\max}(\langle U,U\rangle)}
					$}\\[2pt]
				
				$\mathcal{M}_{\leqslant}^{N;N_{g}}$ &
				Quasi-Stiefel set:
				\par\smallskip
				\centerline{$
					\bigl\{U\in(\mathcal{V}^{N_{g}})^{N}:0<\langle U,U\rangle\leqslant I_{N}\bigr\}
					$}\\[2pt]
				
				$\mathcal{G}_{\leqslant}^{N;N_{g}}$ &
				Quotient set associated with the quasi-Stiefel set:
				\par\smallskip
				\centerline{$
					\mathcal{M}_{\leqslant}^{N;N_{g}}/{\sim}
					$}\\[2pt]
				
				$\operatorname{dist}([\hat{U}],[V_{*}])$ &
				Distance between equivalence classes
				$[\hat{U}],[V_{*}]\in(\mathcal{V}^{N_{g}})^{N}$:
				\par\smallskip
				\centerline{$
					\operatorname{dist}\bigl([\hat{U}],[V_{*}]\bigr)=
					\min_{Q\in\mathcal{O}^{N}}\|\hat{U}Q-V_{*}\|
					$}\\[2pt]
				
				$B(U,\eta)$ &
				Closed $\eta$-neighborhood of
				$U\in(\mathcal{V}^{N_{g}})^{N}$:
				\par\smallskip
				\centerline{$
					B(U,\eta)=\bigl\{\hat{U}\in (\mathcal{V}^{N_{g}})^{N}:\|U-\hat{U}\|\leqslant\eta\bigr\}
					$}\\[2pt]
				
				$B([V_{*}],\eta)$ &
				Closed $\eta$-neighborhood of
				$[V_{*}]\in \mathcal{G}^{N;N_{g}}$:
				\par\smallskip
				\centerline{$
					B([V_{*}],\eta)=\bigl\{\hat{U}\in (\mathcal{V}^{N_{g}})^{N}:
					\operatorname{dist}([\hat{U}],[V_{*}])\leqslant\eta\bigr\}
					$}\\
				\bottomrule[1.2pt]
		\end{tabular}}
	\end{table}

	\section{An iterative method with asymptotic orthogonality}\label{sec:asymptotic-orthogonality algorithm}
	To compute a cluster of eigenpairs of \eqref{linear eigenvalue problem}, we propose a predictor--corrector scheme with asymptotic orthogonality. The {iterates remain in the span of the eigenvectors occurring in the eigenbasis expansions of the initial data}, preserve orthogonality for orthonormal initial data, and achieve asymptotic orthogonality for nonorthogonal initial data: the orthogonality error tends to zero without any orthogonalization operation. {These properties are consistent with those of the continuous model \eqref{eq: continuous model}.}
	
	\subsection{Iterative scheme}
	{Starting from the quasi-Grassmannian gradient flow \eqref{eq: continuous model}, we construct a predictor--corrector time discretization in which the predictor and corrector correspond, respectively, to the terms in the gradient decomposition given below. The resulting time discretization defines an iterative scheme for \eqref{linear eigenvalue problem}.}

	Let $\left\{\tau_n\right\}_{n=0}^{\infty} \subset(0,+\infty)$ be a time-step sequence. We adopt the hypotheses
	\begin{equation*}
		 \sum_{n=0}^{\infty}\tau_n=+\infty.
	\end{equation*}
	The scheme is given by
	\begin{equation}
		\label{Discretization scheme}
		\left\{\begin{aligned}
			&	\hat{U}_{n+1} = U_n -\tau_n \mathcal{A}_{\tilde{U}_{n+\frac{1}{2}}}\tilde{U}_{n+\frac{1}{2}}, \quad \text{where } \tilde{U}_{n+\frac{1}{2}} = \frac{U_n+\hat{U}_{n+1}}{2},
			\\& U_{n+1} = \hat{U}_{n+1} -\tau_n \nabla E(U_{n+1})(I_N-\langle U_{n+1},U_{n+1}\rangle).
		\end{aligned}\right.
	\end{equation}
	{The operator $\mathcal{A}_U$ is defined by
	\begin{equation*}
		\mathcal{A}_UW
		=\nabla E(U)\langle U,W\rangle-U\langle\nabla E(U),W\rangle
		\qquad \forall W\in\left(\mathcal{V}^{N_g}\right)^N.
	\end{equation*}
	The symmetry of the inner product implies that $\mathcal A_U$ is skew-symmetric.} {With this definition, the extended Grassmannian gradient satisfies
	\begin{equation*}
		\nabla_G E(U)
		=
		\mathcal{A}_UU
		+
		\nabla E(U)\left(I_N-\langle U,U\rangle\right).
	\end{equation*}
	The predictor and corrector in \eqref{Discretization scheme} are constructed from the first and second terms in this decomposition, respectively.}

	The scheme ({\ref{Discretization scheme}}) consists of two components:
	\begin{itemize}
		\item \textbf{Predictor step}: This step generates an intermediate approximation $\hat{U}_{n+1}$ by solving the first equation in \eqref{Discretization scheme} while preserving the Gram matrix of \(U_n\) (Lemma \ref{predictor: Gram matrix preservation}).
		\item \textbf{Corrector step}: This correction determines the updated iterate \(U_{n+1}\) from the second equation in \eqref{Discretization scheme} and is used to control the distance to the Stiefel manifold \(\mathcal{M}^{N;N_g}\) (Theorem \ref{orthogonality-error decay}).
		\end{itemize}

		{Throughout this section, we assume that the equations in \eqref{Discretization scheme} admit unique solutions at each iteration. Their existence and uniqueness near \([V_*]\), as needed in the convergence analysis, are established in Lemma \ref{lemma: existence of g}. We denote these solutions by the predictor \(\hat U_{n+1}\) and the corrector \(U_{n+1}\), respectively.}
	
		\begin{figure}
		\centering
		\includegraphics[width=0.59\linewidth]{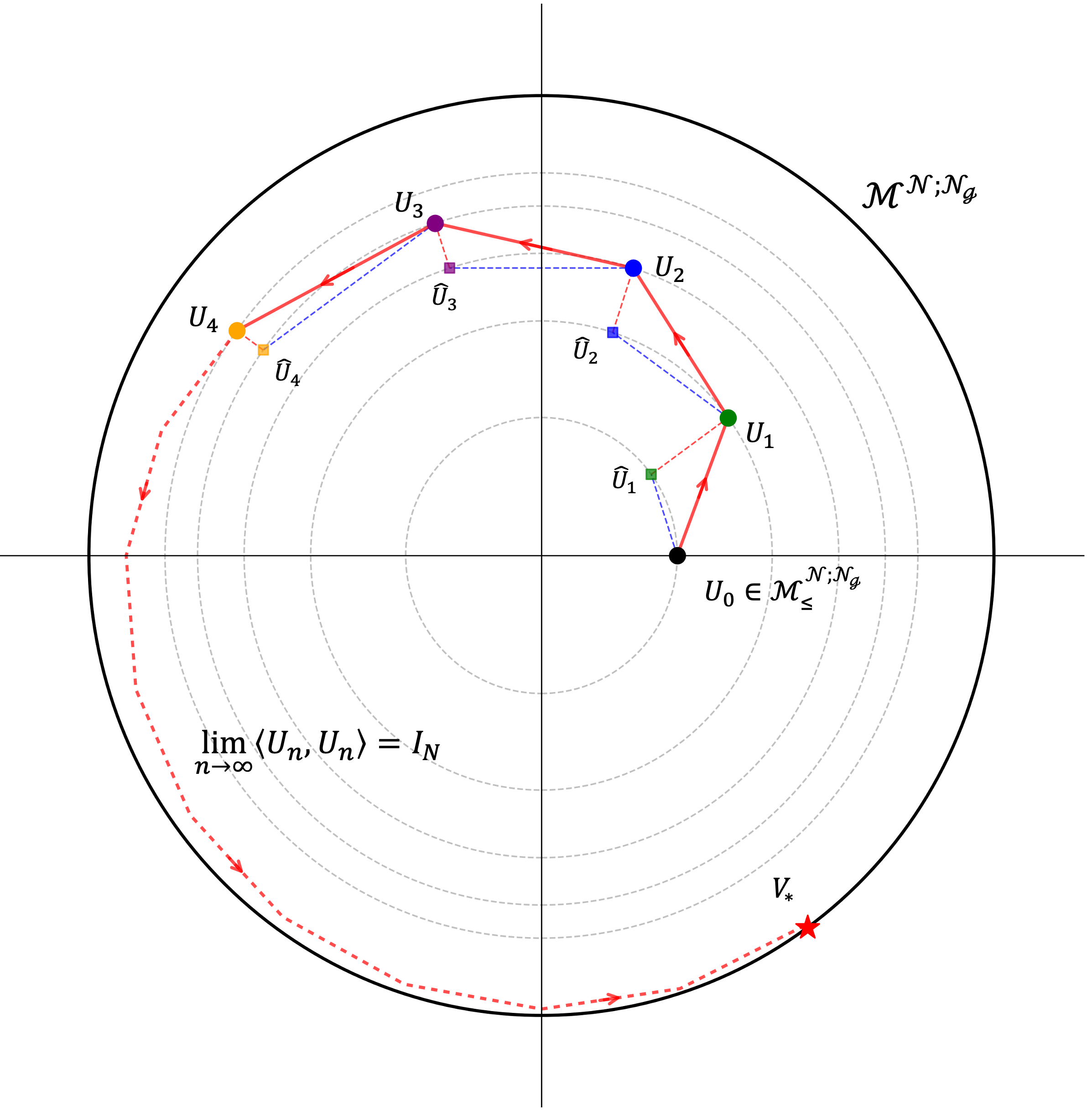}
		\caption{Asymptotic orthogonality of the scheme ({\ref{Discretization scheme}})}
		\label{fig:untendbeortho}
	\end{figure}
	{The analysis below shows that, under appropriate initial conditions, the sequence \(\{U_n\}_{n\geqslant0}\) remains in the quasi-Stiefel set \(\mathcal{M}_{\leqslant}^{N;N_g}\) and satisfies \(\lim\limits_{n\to\infty}\langle U_n,U_n\rangle=I_N\) while approaching the target eigenspace, as illustrated in Fig.~\ref{fig:untendbeortho}.}

	We present the iterative algorithm for \eqref{Discretization scheme} as follows:
	\begin{algorithm}[!h]
		\caption{{Predictor--corrector discretization}}
		\label{alg:Discretization scheme}
		\begin{algorithmic}[1]
			\STATE Given $\epsilon>0$, $\epsilon'>0$, $\delta_T>0$, initial data $U_0\in \left(\mathcal{V}^{N_g}\right)^N$, calculate the gradient and the orthogonality error, and let $n=0$;
			\WHILE{$\left\|\nabla_G E(U_n)\right\|>\epsilon
				\text{ or }
				\left\|I_N-\langle U_n,U_n\rangle\right\|>\epsilon'$}
			\STATE Set time step $\tau_n\leqslant \delta_T$;
			\STATE Solve \eqref{Discretization scheme} to obtain $U_{n+1}$;
			\STATE Let $n=n+1$, and calculate the gradient and the orthogonality error;
			\ENDWHILE
		\end{algorithmic}
\end{algorithm}

{\subsection{Invariant-subspace nonexpansion}}
{Let \(\mathcal V_0\), defined below, be the span of the eigenvectors of \(H\) that occur with nonzero coefficients in the columns of \(U_0\). We prove that the predictor and corrector remain in \((\mathcal V_0)^N\); equivalently, the iteration generates no component orthogonal to \(\mathcal V_0\). We also show that the iterates remain in the quasi-Stiefel set \(\mathcal M_{\leqslant}^{N;N_g}\).}

	The following properties of the predictor follow from Lemmas 4.3 and 4.5 of \cite{dai2020}.
	\begin{lemma}
		\label{predictor: Gram matrix preservation}
		For every \(n=0,1,2,\ldots\), the predictor \(\hat U_{n+1}\) in \eqref{Discretization scheme} satisfies
		\begin{equation*}
				\langle\hat{U}_{n+1},\hat{U}_{n+1}\rangle = \langle U_n,U_n\rangle.
		\end{equation*}
		If, in addition, \(\langle U_n,U_n\rangle\leqslant I_N\), then the spectrum of $\langle\tilde{U}_{n+\frac{1}{2}},\tilde{U}_{n+\frac{1}{2}}\rangle$ satisfies
		\begin{equation*}
				\lambda(\langle\tilde{U}_{n+\frac{1}{2}},\tilde{U}_{n+\frac{1}{2}}\rangle) \subset[0,1],
		\end{equation*}
			where
		\begin{equation}
			\label{eq: Un+1/2}
			\tilde{U}_{n+\frac{1}{2}}  = \left(\mathcal{I} +\frac{\tau_n}{2} \mathcal{A}_{\tilde{U}_{n+\frac{1}{2}}}\right)^{-1} U_n,
			\end{equation}
		and $\mathcal{I}$ denotes the identity operator on \(\left(\mathcal{V}^{N_g}\right)^N\).
		
	\end{lemma}

Let \(v_{\ell_1^0},v_{\ell_2^0},\ldots,v_{\ell_{d_0}^0}\) be the eigenvectors of \(H\) having nonzero coefficients in the eigenbasis expansions of the columns of \(U_0\), and define
\begin{equation*}
	\mathcal{V}_0=\operatorname{span}\{v_{\ell_1^0},v_{\ell_2^0},\cdots,v_{\ell_{d_0}^0}\} .
\end{equation*}
We see that $H\mathcal{V}_0\subset\mathcal{V}_0$, and the notation \(\left(\mathcal{V}_0\right)^N\) denotes the product of this same subspace in all \(N\) columns.

	\begin{lemma}
		\label{lemma: linearly independent+ belong to same space}		
		If $U_n\in   \left(\mathcal{V}_0\right)^N\bigcap \mathcal{M}_{\leqslant}^{N;N_g}$ and $\tau_n<\frac{2}{|\lambda_1-\lambda_{\max}|}$, then
		\begin{equation*}
				\hat{U}_{n+1}\in   \left(\mathcal{V}_0\right)^N \text{ and }U_{n+1}\in   \left(\mathcal{V}_0\right)^N,
		\end{equation*}
		where $\lambda_1<0$ and $\lambda_{\max}$ denote the smallest and largest eigenvalue of $H$, respectively.
	\end{lemma}
	
	\begin{proof}
		We first derive the inner product property of \(\tilde{U}_{n+\frac{1}{2}}\). From (\ref{eq: Un+1/2}) and the skew-symmetry of $\mathcal{A}_{\tilde{U}_{n+\frac{1}{2}}}$,
		\begin{equation*}
			\begin{aligned}
				\langle\tilde{U}_{n+\frac{1}{2}},\tilde{U}_{n+\frac{1}{2}}\rangle & =\left\langle U_n,\left(\mathcal{I} -\frac{\tau_n}{2} \mathcal{A}_{\tilde{U}_{n+\frac{1}{2}}}\right)^{-1}   \left(\mathcal{I} +\frac{\tau_n}{2} \mathcal{A}_{\tilde{U}_{n+\frac{1}{2}}}\right)^{-1} U_n\right\rangle
				\\ &= \left\langle U_n,\left(\mathcal{I} -\frac{\tau_n^2}{4} \mathcal{A}_{\tilde{U}_{n+\frac{1}{2}}}^2\right)^{-1}U_n\right\rangle.
			\end{aligned}
		\end{equation*}
		Since $\mathcal{A}_{\tilde{U}_{n+\frac{1}{2}}}$ is skew-symmetric, $-\mathcal{A}_{\tilde{U}_{n+\frac{1}{2}}}^2$ is positive semidefinite. Therefore
		\begin{equation*}
			\left(\mathcal{I} - \frac{\tau_n^2}{4} \mathcal{A}_{\tilde{U}_{n+\frac{1}{2}}}^2\right)^{-1}>0,
		\end{equation*}
		which together with $\langle U_n,U_n\rangle>0$ leads to $\langle\tilde{U}_{n+\frac{1}{2}},\tilde{U}_{n+\frac{1}{2}}\rangle>0$.
		
		To verify \(\hat{U}_{n+1}, U_{n+1} \in (\mathcal{V}_0)^N\), we proceed by showing that they are orthogonal to any eigenvector outside \(\mathcal{V}_0\). Let \(\tilde{v}\) be an eigenvector of \(H\) with corresponding eigenvalue \(\tilde{\lambda}\), and assume \(\tilde{v}\perp\mathcal{V}_0\). Since \(U_n\in(\mathcal{V}_0)^N\), we have
		\begin{equation*}
			\langle\tilde{v},U_n\rangle=0 .
		\end{equation*}
		Substituting the definition of \(\hat{U}_{n+1}\), we have
		\begin{equation*}
			\begin{aligned}
			&	\langle\tilde{v},\hat{U}_{n+1}\rangle = \langle\tilde{v},U_n\rangle-\tau_n \langle\tilde{v},\mathcal{A}_{\tilde{U}_{n+\frac{1}{2}}}\tilde{U}_{n+\frac{1}{2}}\rangle   
				\\ =&  \langle\tilde{v},U_n\rangle -	\frac{ \tilde{\lambda}\tau_n}{2}\langle\tilde{v},U_n+\hat{U}_{n+1}\rangle\langle\tilde{U}_{n+\frac{1}{2}},\tilde{U}_{n+\frac{1}{2}}\rangle + \frac{\tau_n}{2}\langle\tilde{v},U_n+\hat{U}_{n+1}\rangle\langle\tilde{U}_{n+\frac{1}{2}},H\tilde{U}_{n+\frac{1}{2}}\rangle.
			\end{aligned}
		\end{equation*}
		Rearranging terms yields
		\begin{equation*}
			\begin{aligned}
				&\langle\tilde{v},\hat{U}_{n+1}\rangle \left( I_N + \frac{ \tilde{\lambda}\tau_n }{2}\langle\tilde{U}_{n+\frac{1}{2}},\tilde{U}_{n+\frac{1}{2}}\rangle - \frac{ \tau_n}{2}\langle\tilde{U}_{n+\frac{1}{2}},H\tilde{U}_{n+\frac{1}{2}}\rangle \right)  
				\\ =& \langle\tilde{v},U_n\rangle \left( I_N - \frac{ \tilde{\lambda}\tau_n }{2}\langle\tilde{U}_{n+\frac{1}{2}},\tilde{U}_{n+\frac{1}{2}}\rangle +\frac{ \tau_n}{2}\langle\tilde{U}_{n+\frac{1}{2}},H\tilde{U}_{n+\frac{1}{2}}\rangle \right) . 
			\end{aligned}
		\end{equation*}
		where the right-hand side is zero. For the left-hand side, we use
		\begin{equation*}
			\tau_n<\frac{2}{|\lambda_1-\lambda_{\max}|}
			\leqslant\frac{2}{|\tilde{\lambda}-\lambda_{\max}|}.
		\end{equation*}
		Together with $\langle\tilde{U}_{n+\frac{1}{2}},\tilde{U}_{n+\frac{1}{2}}\rangle>0$, this implies
		\begin{equation*}
			\begin{aligned}
				&	I_N + \frac{ \tilde{\lambda}\tau_n  }{2}\langle\tilde{U}_{n+\frac{1}{2}},\tilde{U}_{n+\frac{1}{2}}\rangle - \frac{ \tau_n}{2}\langle\tilde{U}_{n+\frac{1}{2}},H\tilde{U}_{n+\frac{1}{2}}\rangle 
				\\		\geqslant &\  \left\langle\tilde{U}_{n+\frac{1}{2}},\left( \mathcal{I}+ \frac{ \tilde{\lambda}\tau_n }{2}   \mathcal{I}- \frac{ \tau_n}{2}   H \right)\tilde{U}_{n+\frac{1}{2}}\right\rangle
				\geqslant (1+\frac{\tilde{\lambda}\tau_n }{2}-\frac{ \lambda_{\max} \tau_n}{2} )\langle\tilde{U}_{n+\frac{1}{2}},\tilde{U}_{n+\frac{1}{2}}\rangle>0.
			\end{aligned}
		\end{equation*}
		Thus, $ \left( I_N + \frac{ \tilde{\lambda} \tau_n }{2}\langle\tilde{U}_{n+\frac{1}{2}},\tilde{U}_{n+\frac{1}{2}}\rangle - \frac{ \tau_n}{2}\langle\tilde{U}_{n+\frac{1}{2}},H\tilde{U}_{n+\frac{1}{2}}\rangle \right)  $ is positive definite and hence invertible. This forces \(\langle\tilde{v},\hat{U}_{n+1}\rangle=0\). Since the above reasoning holds for every eigenvector \(\tilde{v}\perp\mathcal{V}_0\), we have
		\begin{equation*}
					\hat{U}_{n+1}\in \left(\mathcal{V}_0\right)^N.
		\end{equation*}
	
	It remains to prove that \(U_{n+1}\in\left(\mathcal{V}_0\right)^N\). Let \(v_k\perp\mathcal{V}_0\) satisfy \(Hv_k=\lambda_kv_k\). Since \(\hat{U}_{n+1}\in\left(\mathcal{V}_0\right)^N\), the corrector equation gives
	\begin{equation*}
		\langle v_k,U_{n+1}\rangle
		\left[
		I_N+\tau_n\lambda_k
		\left(I_N-\langle U_{n+1},U_{n+1}\rangle\right)
		\right]
		=0.
	\end{equation*}
	Set
	\begin{equation*}
		\tilde{U}_{n+1}
		=
		U_{n+1}-2v_k\langle v_k,U_{n+1}\rangle.
	\end{equation*}
	Then
	\begin{equation*}
		\left\langle \tilde{U}_{n+1},\tilde{U}_{n+1}\right\rangle
		=
		\langle U_{n+1},U_{n+1}\rangle
	\end{equation*}
	and
	\begin{equation*}
		\begin{aligned}
			&\tilde{U}_{n+1}
			+\tau_nH\tilde{U}_{n+1}
			\left(
			I_N-\left\langle \tilde{U}_{n+1},\tilde{U}_{n+1}\right\rangle
			\right)
			\\
			={}&
			\hat{U}_{n+1}
			-2v_k\langle v_k,U_{n+1}\rangle
			\left[
			I_N+\tau_n\lambda_k
			\left(I_N-\langle U_{n+1},U_{n+1}\rangle\right)
			\right]
			=\hat{U}_{n+1}.
		\end{aligned}
	\end{equation*}
	The uniqueness of the corrector implies
	\begin{equation*}
		\tilde{U}_{n+1}=U_{n+1},
		\qquad
		\langle v_k,U_{n+1}\rangle=0.
	\end{equation*}
	Therefore,
	\begin{equation*}
		U_{n+1}\in\left(\mathcal{V}_0\right)^N.
	\end{equation*}
	\end{proof}
	
	Let \(\lambda_{\max}^0  = \max\left\{ \frac{(v,Hv)}{(v,v)}:0\ne v\in\mathcal{V}_0 \right\}\). We impose the following assumption for subsequent analysis: \(\lambda_{\max}^0 \leqslant 0\).
	\begin{remark}
		The assumption $\lambda_{\max}^0\leqslant 0 $ can be satisfied. Leveraging the translational invariance of the eigenvalue problem, we can apply a spectral shift to render the relevant eigenvalues of $H$ non-positive; alternatively, we may directly select suitable initial data \(U_0\) to satisfy this condition.
	\end{remark}

	{The following theorem establishes invariant-subspace nonexpansion and preservation of the quasi-Stiefel set.} Under the stated initial-data and time-step conditions, both \(\{U_n\}_{n=0}^{\infty}\) and \(\{\hat{U}_n\}_{n=0}^{\infty}\) remain in
	\begin{equation*}
		\left(\mathcal{V}_0\right)^N
		\bigcap
		\mathcal{M}_{\leqslant}^{N;N_g}.
	\end{equation*}
	\begin{theorem}
		\label{Un inside Stifel manifold}
		Let \(\{U_n\}_{n=0}^\infty\) be the sequence generated by Algorithm \ref{alg:Discretization scheme}. If $\sup\{\tau_n:n\in \mathbb{N} \}\leqslant  \delta_q$ and $U_0 \in \left(\mathcal{V}_0\right)^N\bigcap \mathcal{M}_{\leqslant}^{N;N_g} $, then
		\begin{equation*}
				\hat{U}_n \in  \left(\mathcal{V}_0\right)^N\bigcap \mathcal{M}_{\leqslant}^{N;N_g}  \text{ and }	U_n \in  \left(\mathcal{V}_0\right)^N\bigcap \mathcal{M}_{\leqslant}^{N;N_g}  \quad \forall n\in \mathbb{N}_+,
		\end{equation*}
		where $0<\delta_q <\frac{2}{|\lambda_1-\lambda_{\max}|}$, and $\lambda_1$ and $\lambda_{\max}$ are defined in Lemma \ref{lemma: linearly independent+ belong to same space}.
	\end{theorem}
	\begin{proof}
	We use mathematical induction for the proof. Since $U_0 \in  \left(\mathcal{V}_0\right)^N\bigcap \mathcal{M}_{\leqslant}^{N;N_g} $, we assume $U_n\in  \left(\mathcal{V}_0\right)^N\bigcap \mathcal{M}_{\leqslant}^{N;N_g} $. By Lemma \ref{lemma: linearly independent+ belong to same space}, we first obtain
		\begin{equation*}
				\hat{U}_{n+1},U_{n+1}\in\left(\mathcal{V}_0\right)^N.
		\end{equation*}
		Moreover, Lemma \ref{predictor: Gram matrix preservation} gives
		\begin{equation*}
			\langle\hat{U}_{n+1},\hat{U}_{n+1}\rangle
			=\langle U_n,U_n\rangle>0.
		\end{equation*}
		To show that \(U_{n+1}\) also has full column rank, let \(x\in\mathbb R^N\) satisfy \(U_{n+1}x=0\). Then \(\langle U_{n+1},U_{n+1}\rangle x=0\), and the corrector equation implies
		\begin{equation*}
			\begin{aligned}
				\hat U_{n+1}x
				&=U_{n+1}x+\tau_nHU_{n+1}
				\left(I_N-\langle U_{n+1},U_{n+1}\rangle\right)x\\
				&=0.
			\end{aligned}
		\end{equation*}
		Since \(\langle\hat U_{n+1},\hat U_{n+1}\rangle =\langle U_n,U_n\rangle>0\), this forces \(x=0\). Hence
		\begin{equation*}
			\langle U_{n+1},U_{n+1}\rangle>0.
		\end{equation*}
		The induction hypothesis and Gram matrix preservation also give
		\begin{equation*}
			\langle\hat U_{n+1},\hat U_{n+1}\rangle
			=\langle U_n,U_n\rangle\leqslant I_N,
		\end{equation*}
		and therefore \(\hat U_{n+1}\in(\mathcal V_0)^N \cap\mathcal M_{\leqslant}^{N;N_g}\).
	
		To prove \(\langle U_{n+1},U_{n+1}\rangle\leqslant I_N\), assume to the contrary that \(\langle U_{n+1},U_{n+1}\rangle\) has an eigenvalue \(\mu>1\) with a unit eigenvector \(x\). Then
		$$
		\|U_{n+1}x\|^2=\mu.
		$$
		From the corrector step,
		$$
		\hat{U}_{n+1}x
		=
		U_{n+1}x+\tau_n(1-\mu)HU_{n+1}x.
		$$
		By Lemma \ref{lemma: linearly independent+ belong to same space}, \(U_{n+1}\in(\mathcal{V}_0)^N\). Since \(\lambda_{\max}^0\leqslant0\), the operator \(H\) is non-positive on \(\mathcal{V}_0\). Hence
		$$
		\|\hat{U}_{n+1}x\|
		=
		\|(I+\tau_n(1-\mu)H)U_{n+1}x\|
		\geqslant
		\|U_{n+1}x\|
		=
		\sqrt{\mu}
		>
		1.
		$$
		This contradicts
		$$
		\langle\hat{U}_{n+1},\hat{U}_{n+1}\rangle
		=
		\langle U_n,U_n\rangle
		\leqslant I_N.
		$$
		{The contradiction proves \(\langle U_{n+1},U_{n+1}\rangle\leqslant I_N\), which, together with \(\langle U_{n+1},U_{n+1}\rangle>0\), gives}
		$$
		U_{n+1}\in \left(\mathcal{V}_0\right)^N\bigcap \mathcal{M}_{\leqslant}^{N;N_g}.
		$$
		The proof is complete.
	\end{proof}

	\subsection{Asymptotic orthogonality}
	By Lemma \ref{predictor: Gram matrix preservation} and the uniqueness of the corrector, the scheme preserves orthogonality:
	\begin{equation*}
		U_0\in\mathcal{M}^{N;N_g}
		\quad\Longrightarrow\quad
		U_n\in\mathcal{M}^{N;N_g}
		\qquad \forall n\in\mathbb{N}_+.
	\end{equation*}
	{Even for nonorthogonal initial data, the iterates become asymptotically orthogonal: their Gram matrices converge to $I_N$ without orthogonalization.}
	\begin{theorem}
		\label{orthogonality-error decay}
		Let \(\{U_n\}_{n=0}^\infty\) be the sequence generated by Algorithm \ref{alg:Discretization scheme}. Assume in addition that \(\lambda_{\max}^0<0\). If \(\sup\{\tau_n:n\in \mathbb{N} \} \leqslant \delta_I\) and $U_0 \in \left(\mathcal{V}_0\right)^N\bigcap \mathcal{M}_{\leqslant}^{N;N_g} $, then
		\begin{equation*}
			\begin{aligned}
				\left\|I_N -\langle U_n,U_n\rangle\right\|  
				\leqslant \left(\prod_{k=0}^{n-1}\frac{1}{1+\tau_k|\lambda_{\max}^0|\lambda_{\min}(\langle U_0,U_0\rangle)}\right)
				\|I_N - \langle U_0,U_0\rangle\| \qquad \forall n\in \mathbb{N}_{+}.
			\end{aligned}
		\end{equation*}
		Consequently,
		\begin{equation*}
			\label{discretization tends to I} 
				\lim\limits_{n\rightarrow \infty } \langle U_n,U_n\rangle = I_N. 
		\end{equation*}
		Here $0<\delta_{I} <\min\left\{\delta_q,\frac{1}{|\lambda_1|}\right\}$.
	\end{theorem}
	\begin{proof}
		By Theorem \ref{Un inside Stifel manifold}, \(U_n\in \mathcal{M}_{\leqslant}^{N;N_g}\) for all \(n\in\mathbb{N}\). Thus \(0<\langle U_{n+1},U_{n+1}\rangle\leqslant I_N\). Let \(x_i\) be a unit eigenvector of \(\langle U_{n+1},U_{n+1}\rangle\) with eigenvalue \(\mu_i\). Then \(0<\mu_i\leqslant 1\), and the corrector equation gives
		\begin{equation*}
			\begin{aligned}
				\hat{U}_{n+1}x_i
				&=U_{n+1}x_i+\tau_nHU_{n+1}(I_N-\langle U_{n+1},U_{n+1}\rangle)x_i
				\\&=\left(I+\tau_n(1-\mu_i)H\right)U_{n+1}x_i .
			\end{aligned}
		\end{equation*}
		Since \(U_{n+1}\in \left(\mathcal{V}_0\right)^N\), the vector \(U_{n+1}x_i\) belongs to \(\mathcal{V}_0\). On \(\mathcal{V}_0\), the eigenvalues of \(H\) lie in \([\lambda_1,\lambda_{\max}^0]\). Since \(\lambda_{\max}^0<0\) and \(\tau_n\leqslant \delta_I<1/|\lambda_1|\), we have
		\begin{equation*}
			\left\|\left(I+\tau_n(1-\mu_i)H\right)v\right\|
			\leqslant \left(1-\tau_n(1-\mu_i)|\lambda_{\max}^0|\right)\|v\|
			\qquad \forall v\in \mathcal{V}_0 .
		\end{equation*}
		Using Lemma \ref{predictor: Gram matrix preservation}, \(\langle\hat{U}_{n+1},\hat{U}_{n+1}\rangle=\langle U_n,U_n\rangle\). Moreover,
		\begin{equation*}
			0<1-\tau_n(1-\mu_i)|\lambda_{\max}^0|\leqslant1,
		\end{equation*}
		Consequently,
		\begin{equation}\label{ieq:estimate of mu}
			\begin{aligned}
				x_i^\top\langle U_n,U_n\rangle x_i
				&=\|\hat{U}_{n+1}x_i\|^2
				\\&\leqslant \left(1-\tau_n(1-\mu_i)|\lambda_{\max}^0|\right)^2\|U_{n+1}x_i\|^2
				\\&\leqslant \left(1-\tau_n(1-\mu_i)|\lambda_{\max}^0|\right)\mu_i .
			\end{aligned}
		\end{equation}
		It follows that
		\begin{equation*}
			\begin{aligned}
				x_i^\top(I_N-\langle U_n,U_n\rangle)x_i
				&\geqslant 1-\left(1-\tau_n(1-\mu_i)|\lambda_{\max}^0|\right)\mu_i
				\\&=(1-\mu_i)\left(1+\tau_n|\lambda_{\max}^0|\mu_i\right).
			\end{aligned}
		\end{equation*}
		Since \(1-\tau_n(1-\mu_i)|\lambda_{\max}^0|\leqslant1\), the same estimate \eqref{ieq:estimate of mu} gives
		\begin{equation*}
			x_i^\top\langle U_n,U_n\rangle x_i\leqslant\mu_i.
		\end{equation*}
		Therefore
		\begin{equation*}
			\lambda_{\min}(\langle U_{n+1},U_{n+1}\rangle)\geqslant \lambda_{\min}(\langle U_n,U_n\rangle),
		\end{equation*}
		and by induction,
		\begin{equation*}
			\mu_i\geqslant \lambda_{\min}(\langle U_0,U_0\rangle).
		\end{equation*}
		Consequently,
		\begin{equation*}
			\begin{aligned}
				1-\mu_i
				&\leqslant
				\frac{x_i^\top(I_N-\langle U_n,U_n\rangle)x_i}
				{1+\tau_n|\lambda_{\max}^0|\lambda_{\min}(\langle U_0,U_0\rangle)}
				\\&\leqslant
				\frac{\|I_N-\langle U_n,U_n\rangle\|}
				{1+\tau_n|\lambda_{\max}^0|\lambda_{\min}(\langle U_0,U_0\rangle)} .
			\end{aligned}
		\end{equation*}
		Taking the maximum over all \(i\), we obtain
		\begin{equation*}
			\|I_N-\langle U_{n+1},U_{n+1}\rangle\|
			\leqslant
			\frac{1}{1+\tau_n|\lambda_{\max}^0|\lambda_{\min}(\langle U_0,U_0\rangle)}
			\|I_N-\langle U_n,U_n\rangle\|.
		\end{equation*}
		Iterating this inequality gives
		\begin{equation*}
			\left\|I_N -\langle U_n,U_n\rangle\right\|
			\leqslant \left(\prod_{k=0}^{n-1}\frac{1}{1+\tau_k|\lambda_{\max}^0|\lambda_{\min}(\langle U_0,U_0\rangle)}\right)
			\|I_N - \langle U_0,U_0\rangle\|.
		\end{equation*}
		Since \(U_0\in \mathcal{M}_{\leqslant}^{N;N_g}\), \(\lambda_{\min}(\langle U_0,U_0\rangle)>0\). The time-step condition \(\sum_{k=0}^{\infty}\tau_k=+\infty\) then implies
		\begin{equation*}
			\lim\limits_{n\rightarrow \infty } \langle U_n,U_n\rangle = I_N .
		\end{equation*}
		The proof is complete.
	\end{proof}

	{\begin{corollary}
		\label{corollary: exponential orthogonality decay}
		Under the assumptions of Theorem \ref{orthogonality-error decay}, if \(\tau_n\geqslant \tau_*>0\) for all \(n\in\mathbb N\), then
		\begin{equation*}
			\left\|I_N-\langle U_n,U_n\rangle\right\|
			\leqslant
			\omega^n\left\|I_N-\langle U_0,U_0\rangle\right\|,
		\end{equation*}
		where
		\begin{equation*}
			\omega=
			\left(1+\tau_*|\lambda_{\max}^0|
			\lambda_{\min}(\langle U_0,U_0\rangle)\right)^{-1}\in(0,1).
		\end{equation*}
	\end{corollary}

	Theorem \ref{orthogonality-error decay} establishes asymptotic orthogonality, while Corollary \ref{corollary: exponential orthogonality decay} gives an exponential decay estimate for the orthogonality error.} Moreover, if \(U_k\in\mathcal{M}^{N;N_g}\) for some \(k\in\mathbb{N}\), then
	\begin{equation*}
		U_n\in\mathcal{M}^{N;N_g}
		\qquad \forall n\geqslant k.
	\end{equation*}

	\begin{remark}
		{The scheme in \cite{dai2020} preserves the current Gram matrix, whereas, for the nonorthogonal initial data covered by Theorem \ref{orthogonality-error decay}, the present iteration reduces the orthogonality error.}
	\end{remark}

	\section{Convergence}\label{sec:Numerical analysis}
	We first state four auxiliary results and then prove convergence of Algorithm \ref{alg:Discretization scheme} from a neighborhood of the target eigenspace. {The proofs of the following four auxiliary lemmas are collected in Appendix \ref{app:auxiliary-proofs}.}

	Following \cite[Lemma 4.2]{dai2020}, the first result uses the implicit function theorem to establish {existence and uniqueness of solutions to \eqref{Discretization scheme}} near \([V_*]\).
	\begin{lemma}\label{lemma: existence of g}
		There exist constants $0<\eta_a<\eta_b$ and $\delta^*>0$ and unique continuously differentiable mappings
		\begin{equation*}
			\hat{g},g:B\left([V_*],\eta_a\right)\times[0,\delta^*]
			\longrightarrow B\left([V_*],\eta_b\right)
		\end{equation*}
		such that $\hat{g}(U,0)=g(U,0)=U$ and
		\begin{equation*}
			\begin{aligned}
				&	\hat{g}(U,\tau)
				=U-\tau\mathcal{A}_{\frac{\hat{g}(U,\tau)+U}{2}}\frac{\hat{g}(U,\tau)+U}{2},
				\\&		g(U,\tau)=\hat{g}(U,\tau)-\tau\nabla E(g(U,\tau))
				\left(I_N-\langle g(U,\tau),g(U,\tau)\rangle\right).
			\end{aligned}
		\end{equation*}
		For each $U\in B([V_*],\eta_a)$ and $\tau\in[0,\delta^*]$, $\hat{g}(U,\tau)$ is the unique solution of the predictor equation in a neighborhood of $U$, and $g(U,\tau)$ is the unique solution of the corrector equation in a neighborhood of $\hat{g}(U,\tau)$. If $U\in B(V_*,\eta_a)$, then $\hat{g}(U,\tau),g(U,\tau)\in B(V_*,\eta_b)$. Moreover,
		\begin{equation*}
			\hat{g}(UQ,\tau)=\hat{g}(U,\tau)Q,
			\qquad
			g(UQ,\tau)=g(U,\tau)Q,
		\end{equation*}
		for every $Q\in\mathcal{O}^N$.
	\end{lemma}
	
	The mappings in Lemma \ref{lemma: existence of g} allow us to write the update in \eqref{Discretization scheme} as
	\begin{equation*}
		\hat{U}_{n+1}=\hat{g}(U_n,\tau_n),
		\qquad
		U_{n+1}=g(U_n,\tau_n).
	\end{equation*}
	We next establish the energy estimate used in the convergence analysis.
	\begin{lemma}
		\label{lemma of energy decreases}
		There holds
		\begin{equation*}
			\begin{aligned}
				E(U)-	E\left(g(U, \tau)\right) &\geqslant \tau \left(\frac{1}{2} - \frac{ \tau}{2}|\lambda_1|\right) \left\| \mathcal{A}_{ \frac{\hat{g}(U , \tau)+U }{2}} \frac{\hat{g}(U , \tau)+U }{2} \right\|^2,
				\\&\qquad\qquad\qquad \forall U \in B\left(V_{*}, \eta_a\right)\bigcap \left(\mathcal{V}_0\right)^N\bigcap \mathcal{M}_{\leqslant}^{N;N_g} ,\ \tau \in [0, \delta_e],
			\end{aligned}
		\end{equation*}
			where \(0<\delta_e<\min\{ \delta_q, \delta^*,1/|\lambda_1|\}\), with $\delta_q$ and $\delta^* $ defined in Theorem \ref{Un inside Stifel manifold} and Lemma \ref{lemma: existence of g}, respectively.
	\end{lemma}

	To control the displacement of the predictor in a neighborhood of \([V_*]\), we use the following Lipschitz estimates for \(\mathcal{A}_U\).
	\begin{lemma}
		\label{extended gradient is Lip}
	There holds
		\begin{align}
			\| \mathcal{A}_{U_i}U_i - 	\mathcal{A}_{U_j}U_j\| \leqslant L \|U_i - U_j\|, \quad \forall U_i, U_j \in  B\left(V_{*}, \max\{\eta_a, \eta_b \} \right) \nonumber,
			\\	\|\mathcal{A}_{U_i} - \mathcal{A}_{U_j}\| \leqslant \hat{L}\|U_i - U_j\|, \quad \forall U_i, U_j \in  B\left(V_{*}, \max\{\eta_a, \eta_b \} \right),	\label{Lip of mathcal_A}
		\end{align}
		where $L = 6\alpha^2 \scalebox{0.9}{$ \max\{|\lambda_1 | , | \lambda_{\max} | \} $} $, $\hat{L} =  \frac{2L}{3\alpha}$, and $\alpha = \scalebox{0.85}{$ \max\left\{\sigma_{\max}(U): U \in B\left(V_{*}, 	\scalebox{0.9}{$ \max\{\eta_a, \eta_b \} $} \right)\right\}  $}$.
	\end{lemma}

	Combining the preceding well-posedness and Lipschitz estimates gives the bound required to keep the iteration in this neighborhood.
	\begin{lemma}
		\label{discretization is bounded}
		There holds
		\begin{equation*} 
				g\big(B([V_{*}],\eta_2)\times [0, \delta_b] \big)\subset B([V_{*}],\eta_1),
		\end{equation*}
		where $\eta_1 \in \left(0, \min\{1, \eta_a, \eta_b\}\right)$, $\eta_2\in (0,\eta_1)$, and
		\begin{equation*}
			\delta_b\leqslant \min\left\{
			\frac{\eta_1-\eta_2}
			{L\max\{\eta_a,\eta_b\}+|\lambda_1|\alpha(\alpha+1)\eta_b},
			\delta^*\right\}.
		\end{equation*}
	\end{lemma}

	Note that $[V_*]$ is the unique minimizer of $E$ in
	$$
	B\left(\left[V_{*}\right], \eta_a\right)\bigcap\mathcal{G}^{N;N_g}.
	$$
	Under $\lambda_N<0$, $[V_*]$ is also the unique minimizer in
	$$
	B\left(\left[V_{*}\right], \eta_a\right) \bigcap \left(\mathcal{V}_0\right)^N \bigcap\mathcal{G}^{N;N_g}_{\leqslant}.
	$$
	To facilitate subsequent analysis, we define
	\begin{equation*}
		E_0= \min\left\{ E(\tilde{U}): [\tilde{U}] \subset \overline{ B([V_{*}], \eta_1)\backslash B([V_{*}], \eta_2)}\bigcap \left(\mathcal{V}_0\right)^N\bigcap\mathcal{G}^{N;N_g}_{\leqslant}\right\},
	\end{equation*}
	for some $\eta_1 \in \left(0, \min\{1, \eta_a, \eta_b\}\right)$ and $\eta_2\in (0,\eta_1)$. We further define
	\begin{equation*} 
			\mathcal{S} = \left\{U \in (\mathcal{V}^{N_g})^N: [U] \subset B([V_{*}], \eta_2)\bigcap \left(\mathcal{V}_0\right)^N \bigcap\mathcal{G}^{N;N_g}_{\leqslant}\bigcap \mathcal{L}_{E_1}\right\}, 
	\end{equation*}
	where $\mathcal{L}_{E_1} =  \left\{U \in (\mathcal{V}^{N_g})^N: E(U) \leqslant \frac{ E(V_{*})+ E_0}{2} \stackrel{\Delta}{=}E_1\right\}$.
	
	We state the convergence of Algorithm \ref{alg:Discretization scheme} as follows:
	\begin{theorem}
		\label{main theorem of convergence}
		{Assume that \(\lambda_{\max}^0<0\) and \(U_0\in\mathcal S\). Set
		\begin{equation*}
			\delta_T=\min\{\delta_I,\delta_e,\delta_b\},
		\end{equation*}
		where \(\delta_I\), \(\delta_e\), and \(\delta_b\) are defined in Theorem \ref{orthogonality-error decay}, Lemma \ref{lemma of energy decreases}, and Lemma \ref{discretization is bounded}, respectively. Then the sequence \(\{U_n\}_{n=0}^{\infty}\) generated by Algorithm \ref{alg:Discretization scheme} satisfies}
		\begin{equation}
				\label{discretization energy decrease} 
				E(U_n)-E(U_{n+1})  \geqslant \tau_n\left(\frac{1}{2} - \frac{\tau_n}{2}|\lambda_1|\right)\|\mathcal{A}_{\tilde{U}_{n+\frac{1}{2}}}\tilde{U}_{n+\frac{1}{2}}\|^2 , 
		\end{equation}
		along with the asymptotic behaviors:
		\begin{equation*}
			\begin{aligned}
				&\lim\limits_{n\rightarrow\infty} \left\|\nabla_G E(U_n)\right\| = 0,
				\\&\lim\limits_{n\rightarrow\infty}E(U_n) = E(V_{*}),
				\\&\lim_{n \rightarrow\infty}\operatorname{dist}\left([U_n],[V_{*}]\right)=0.
			\end{aligned}
		\end{equation*}
	\end{theorem}
	\begin{proof}
		We first prove that \(U_n\in \mathcal{S}\) for all \(n\in\mathbb{N}\). This holds for \(n=0\) by assumption. Suppose \(U_n\in\mathcal{S}\). Since \(\delta_T\leqslant\delta_e\), Lemma \ref{lemma of energy decreases} gives
		\begin{equation*}
				E(U_n)-E(U_{n+1})  \geqslant \tau_n\left(\frac{1}{2} - \frac{\tau_n}{2}|\lambda_1|\right)\|\mathcal{A}_{\tilde{U}_{n+\frac{1}{2}}}\tilde{U}_{n+\frac{1}{2}}\|^2 .
		\end{equation*}
		Thus \(E(U_{n+1})\leqslant E(U_n)\leqslant E_1\). Since \(\delta_T\leqslant\delta_b\), Lemma \ref{discretization is bounded} gives
		\begin{equation*}
				U_{n+1}\in B\left(\left[V_{*}\right], \eta_1\right).
		\end{equation*}
		Moreover, Theorem \ref{Un inside Stifel manifold} implies
		\begin{equation*}
				[U_{n+1}]\subset \left(\mathcal{V}_0\right)^N\bigcap \mathcal{G}^{N;N_g}_{\leqslant}.
		\end{equation*}
		If \(U_{n+1}\notin B([V_*],\eta_2)\), then the definition of \(E_0\) gives \(E(U_{n+1})\geqslant E_0\), which contradicts \(E(U_{n+1})\leqslant E_1<E_0\). Hence \(U_{n+1}\in\mathcal{S}\). By induction, \(U_n\in\mathcal{S}\) for all \(n\in\mathbb{N}\), and the previous estimate gives \eqref{discretization energy decrease}.

		Since \(\tau_n\leqslant\delta_e<1/|\lambda_1|\), the factor \(\frac{1}{2}-\frac{\tau_n}{2}|\lambda_1|\) is bounded from below by \(\frac{1}{2}-\frac{\delta_e}{2}|\lambda_1|>0\). Using the lower boundedness of \(\{E(U_n)\}\), we get
		\begin{equation*} 
				\left(\frac{1}{2} - \frac{\delta_e}{2}|\lambda_1|\right)
				\sum\limits_{n=0}^{\infty} \tau_n
				\left\|\mathcal{A}_{\tilde{U}_{n+\frac{1}{2}}}\tilde{U}_{n+\frac{1}{2}}\right\|^2
				\leqslant E(U_0)-\lim\limits_{n\rightarrow \infty }E(U_n), 
		\end{equation*}
		which together with $\sum\limits_{n=0}^{\infty}  \tau_n = \infty$ implies
		\begin{equation*} 
				\liminf\limits_{n\rightarrow\infty}	\left\|\mathcal{A}_{\tilde{U}_{n+\frac{1}{2}}}\tilde{U}_{n+\frac{1}{2}}\right\| = 0. 
		\end{equation*}
		Consequently, due to (\ref{Discretization scheme}), there exists a subsequence $\{\tilde{U}_{n_{k}+\frac{1}{2}}\}_{k=0}^{\infty}$, such that
		\begin{equation*} 
				\lim\limits_{k\rightarrow\infty}\left\|\hat{U}_{n_{k}+1} - U_{n_k}\right\|\leqslant \delta_T\lim\limits_{k\rightarrow\infty} 	\left\|\mathcal{A}_{\tilde{U}_{n_{k}+\frac{1}{2}}}\tilde{U}_{n_{k}+\frac{1}{2}}\right\|  = 0. 
		\end{equation*}
		We obtain from Theorem \ref{orthogonality-error decay} and Lemma \ref{predictor: Gram matrix preservation} that
		\begin{equation*} 
				\lim\limits_{k\rightarrow\infty}\left\|I_N - \langle\hat{U}_{n_k+1},\hat{U}_{n_k+1}\rangle\right\|=	\lim\limits_{k\rightarrow\infty}\left\|I_N - \langle U_{n_k},U_{n_k}\rangle\right\| = 0. 
		\end{equation*}
		Using the triangle inequality and \eqref{Discretization scheme} again:
		\begin{equation*}
			\begin{aligned}
				\lim\limits_{k\rightarrow\infty}\left\|U_{n_{k}+1} - U_{n_k}\right\|
				\leqslant&  \lim\limits_{k\rightarrow\infty}\left\|U_{n_{k}+1} -\hat{U}_{n_{k}+1} \right\|+	\lim\limits_{k\rightarrow\infty}\left\|\hat{U}_{n_{k}+1} - U_{n_k}\right\|
				\\	\leqslant & \delta_T\sup_{\tilde{U}\in\mathcal{S}}\left\|\nabla E(\tilde{U}) \right\|\lim\limits_{k\rightarrow\infty}\left\| I_N - \langle U_{n_{k}+1},U_{n_{k}+1}\rangle\right\|+0  = 0.
			\end{aligned}
		\end{equation*}
		Noting that $\mathcal{S}$ is compact, we have a subsequence of $\{U_{n_k}\}_{k=0}^{\infty}$, which is still denoted as $\{U_{n_k}\}_{k=0}^{\infty}$ for simplicity, satisfying
		\begin{equation*}
			\begin{aligned}
				\lim_{k \rightarrow \infty}U_{n_k} = \bar{U}
			\end{aligned}
		\end{equation*}
		for some $\bar{U} \in \mathcal{S}$. It follows that
		\begin{equation*} 
				\lim_{k \rightarrow \infty} \hat{U}_{n_k+1} =  \bar{U}, 
		\end{equation*}
		and
		\begin{equation*} 
				\lim_{k \rightarrow \infty} \tilde{U}_{n_{k}+\frac{1}{2}} =\lim_{k \rightarrow \infty} \frac{\hat{U}_{n_k+1}+U_{n_k}}{2} =  \bar{U}. 
		\end{equation*}
		Taking limits in \(\mathcal{A}_{\tilde{U}_{n_{k}+\frac{1}{2}}} \tilde{U}_{n_{k}+\frac{1}{2}} \to 0\), we get
		\begin{equation*} 
				\nabla_{G}E(\bar{U}) =\mathcal{A}_{\bar{U}}\bar{U} = 0, 
		\end{equation*}
		which implies \([\bar{U}]\) is a critical point of $E$. 
		Because of the uniqueness of the critical point in $B([V_{*}], \eta_1)$, we have $[\bar{U}] = [V_{*}]$ and
		\begin{equation*} 
				\lim_{n \rightarrow\infty}E(U_n) = \lim_{k \rightarrow \infty}E(U_{n_k}) = E(\bar{U}). 
		\end{equation*}
		
		To prove \(\lim_{n \to \infty} \operatorname{dist}([U_n], [V_{*}]) = 0\), suppose for contradiction there exists a subsequence $\{U_{n_l}\}_{l=0}^{\infty}$ and $\hat{\delta}>0$ such that $\operatorname{dist}([U_{n_l}],[V_{*}])>\hat{\delta}$.  By compactness of $\mathcal{S}$, we have a subsequence of $\{U_{n_l}\}_{l=0}^{\infty}$, which is still denoted as $\{U_{n_l}\}_{l=0}^{\infty}$ for simplicity, satisfying $	\lim\limits_{l\rightarrow \infty}U_{n_l} =\check{U}$ for some $\check{U} \in \mathcal{S}$. Thus we have
		\begin{equation*} 
				E(\check{U}) = \lim\limits_{l \rightarrow \infty}E(U_{n_l})  = E(V_{*}). 
		\end{equation*}
		By  the uniqueness of the minimizer in $B\left([V_{*}],  \eta_1\right)$, we obtain $[\check{U}] = [V_{*}]$, which contradicts the assumption $\operatorname{dist}\left([U_{n_l}],[V_{*}]\right)>\hat{\delta}$.
		
		Finally, since \(\nabla_G E\) is continuous on \(\mathcal{S}\) and \(\nabla_G E(V_{*}) = 0\), we obtain
		\begin{equation*} 
				\lim\limits_{n\rightarrow \infty}\left\|\nabla_{G}E(U_n) \right\|  = \left\|\nabla_{G}E(V_{*}) \right\| =0. 
		\end{equation*}
		
		This completes the proof.
	\end{proof}

		\section{Numerical experiments}\label{Numerical experiments}
		{For the numerical realization of \eqref{Discretization scheme},we perform $p_n$ fixed-point iterations for the predictor and perform a single Picard iteration for the corrector. The resulting procedure is summarized in Algorithm \ref{alg:Practical Iteration}, and its derivation and implementation details are given in Appendix \ref{app:practical-implementation}.} We report three numerical experiments whose results exhibit energy decrease, convergence behavior consistent with the analysis in Sections \ref{sec:asymptotic-orthogonality algorithm} and \ref{sec:Numerical analysis}, and decreasing orthogonality errors.

		\begin{algorithm}[H]
			\caption{Practical iteration}
			\label{alg:Practical Iteration}
			\begin{algorithmic}[1]
				\STATE Given $\epsilon>0$, $\epsilon'>0$, $\tilde{\delta}_T>0$, initial data $U_0\in (\mathcal{V}^{N_g})^N$, calculate the gradient and the orthogonality error, and let $n=0$;
				\WHILE{$\left\|\nabla_G E(U_n)\right\|>\epsilon
					\text{ or }
					\left\|I_N-\langle U_n,U_n\rangle\right\|>\epsilon'$}
				\STATE Set time step $\tau_n\leqslant \tilde{\delta}_T$ and choose $p_n \in \mathbb{N}_{+}$;
				\STATE $\tilde{U}_{n+\frac{1}{2}}^{(0)}=U_n$;
				\FOR{$k=1, \ldots, p_n$}
				\STATE
			 \qquad $	\tilde{U}_{n+\frac{1}{2}}^{(k)} = \left(\mathcal{I} +\frac{\tau_n}{2}\mathcal{A}_{\tilde{U}_{n+\frac{1}{2}}^{(k-1)}}\right)^{-1} U_n; $
				\ENDFOR
				\STATE $\hat{U}_{n+1}  = 2\tilde{U}_{n+\frac{1}{2}}^{(p_n)} -U_n$;
				\STATE Perform one Picard iteration for the corrector with initial value \(\hat{U}_{n+1}\):
				\STATE \qquad $U_{n+1} = \hat{U}_{n+1} -\tau_n\nabla E(\hat{U}_{n+1})(I_N - \langle U_n,U_n\rangle)$;
				\STATE Let $n=n+1$, and calculate the gradient and the orthogonality error;
				\ENDWHILE
			\end{algorithmic}
		\end{algorithm}

	We test Algorithm \ref{alg:Practical Iteration} on three eigenvalue problems governed by the Laplacian, harmonic oscillator, and hydrogen atom Schr\"odinger operators. All experiments are conducted on the LSSC-IV platform at the Academy of Mathematics and Systems Science, Chinese Academy of Sciences. Spatial discretization is implemented using the finite element method with quadratic elements. {We use two fixed-point updates for the predictor in each iteration, namely, $p_n=2$.}

	Motivated by \cite{dai2019adaptive,dai2021convergent}, we use the following adaptive time step:
	\begin{equation*} 
			\tau_n=\min\left\{	\frac{\left\|\nabla_G E\left(U_n\right)\right\|^2}{\operatorname{Hess}_{\tilde{G}}\left(U_n\right)\left[\nabla_G E\left(U_n\right), \nabla_G E\left(U_n\right)\right]},\tilde{\delta}_T \right\}, \qquad n=0,1,2,\ldots, 
	\end{equation*}
	where, for all $U,V,W\in\big(\mathcal{V}^{N_g}\big)^N$, the Hessian extension to \((\mathcal{V}^{N_g})^N\) \cite{Wang2025} is defined by
	\begin{equation*}
		\operatorname{Hess}_{\tilde{G}}(U)[V,W]
		=\left(\nabla^2 E(U)W-W\langle U,\nabla E(U)\rangle,V\right).
	\end{equation*}
	This strategy remains effective even if \(U_n \notin \mathcal{M}^{N;N_g}\).

	Reference solutions \((V_*, \Lambda_*)\) are computed via the \texttt{eigs} solver from \texttt{Arpack.jl} for (\ref{linear eigenvalue problem}). The initial iterate \(U_0\) is randomly generated with linearly independent columns.\footnote{The entries of the finite element coefficient vector of each column are sampled independently from the uniform distribution on \([-\frac{1}{2},\frac{1}{2})\). No condition \(\langle U_0,U_0\rangle\leqslant I_N\) or mutual orthogonality is imposed. These initial data are therefore less restrictive than those covered by the analysis; extending the convergence theory to this setting remains open.} Unless stated otherwise, the iteration terminates only when
	\begin{equation*}
		\left\|\nabla_G E(U_n)\right\|<10^{-5}
		\qquad\text{and}\qquad
		\left\|I_N-\langle U_n,U_n\rangle\right\|<10^{-10}.
	\end{equation*}
	{That is, we choose \(\epsilon=10^{-5}\) and \(\epsilon'=10^{-10}\) in these experiments.}
	
	For clarity, we define the relative difference from the final computed iterate by
	\begin{equation}\label{eq:relative error of iterates} 
			\text{err}_{U_n}=\frac{\|U_{n}-U_{\text{end}}\|}{\|U_{\text{end}}\|}, 
	\end{equation}
	where \(U_{\text{end}}\) denotes the final iteration result. This quantity measures the stabilization of \(U_n\) relative to \(U_{\mathrm{end}}\).
\begin{remark}
		Since $U_{\mathrm{end}}$ is numerically orthogonal, let $Q_{\mathrm{end}}\in\mathcal{O}^N$ and \(\Lambda_{\mathrm{end}}=\operatorname{diag} (\lambda_{1,\mathrm{end}},\ldots,\lambda_{N,\mathrm{end}})\) satisfy
	\begin{equation*}
		\left\langle U_{\mathrm{end}},HU_{\mathrm{end}}\right\rangle
		Q_{\mathrm{end}}=Q_{\mathrm{end}}\Lambda_{\mathrm{end}}.
	\end{equation*}
	The eigenvalue approximations are the diagonal entries of \(\Lambda_{\mathrm{end}}\), and the corresponding eigenvector approximations are \(V_{\mathrm{end}}=U_{\mathrm{end}}Q_{\mathrm{end}}\). Since $Q_{\mathrm{end}}\in\mathcal{O}^N$, \([V_{\mathrm{end}}]=[U_{\mathrm{end}}]\), so $V_{\mathrm{end}}$ is a representative of the computed approximation to the target eigenspace \([V_*]\). Its residual satisfies
	\begin{equation*}
		\begin{aligned}
			HV_{\mathrm{end}}-V_{\mathrm{end}}\Lambda_{\mathrm{end}}
			&=\left(HU_{\mathrm{end}}-U_{\mathrm{end}}
			\left\langle U_{\mathrm{end}},HU_{\mathrm{end}}\right\rangle\right)
			Q_{\mathrm{end}}
			\\&=\nabla_G E(U_{\mathrm{end}})Q_{\mathrm{end}},
		\end{aligned}
	\end{equation*}
	and hence
	\begin{equation*}
		\left\|HV_{\mathrm{end}}
		-V_{\mathrm{end}}\Lambda_{\mathrm{end}}\right\|
		=\left\|\nabla_G E(U_{\mathrm{end}})\right\|.
	\end{equation*}
	Thus, the gradient norm is the residual norm of the corresponding eigenpair approximations.
\end{remark}

		\begin{example}\label{eq:3D Laplace eigenvalue equation}
			We consider the Laplacian eigenvalue problem on the bounded domain \(\Omega = (0, \pi)^3\) with homogeneous Dirichlet boundary conditions \cite{cheng1975eigenfunctions}:
			\begin{equation*}
				-\frac{1}{2} \Delta u = \lambda u \quad \text{in} \ \Omega, \qquad \int_{\Omega} u^2 = 1,
			\end{equation*}
			where \((u, \lambda) \in H_0^1(\Omega) \times \mathbb{R}\). The explicit eigenvalues and eigenfunctions are
			\begin{equation*} 
					\lambda_{k_1, k_2, k_3} = \frac{k_1^2 + k_2^2 + k_3^2}{2}, 
					\quad u_{k_1, k_2, k_3}(x) = \left( \frac{2}{\pi} \right)^{\frac{3}{2}} \sin(k_1 x_1) \sin(k_2 x_2) \sin(k_3 x_3), 
			\end{equation*}
			for $k_1, k_2, k_3\in \mathbb{N}_+$. {We choose \(N=11\) to include the first five distinct eigenvalues with their full multiplicities, so that \(\lambda_N<\lambda_{N+1}\).} We use a uniform finite element mesh with \(N_g = 24389\) degrees of freedom. Numerical results are shown in Fig.~\ref{fig:laplace_total}. The main numerical observations are summarized below:
			\begin{figure}[htbp]
				\centering
				\begin{subfigure}[b]{0.48\linewidth}
					\includegraphics[width=\linewidth]{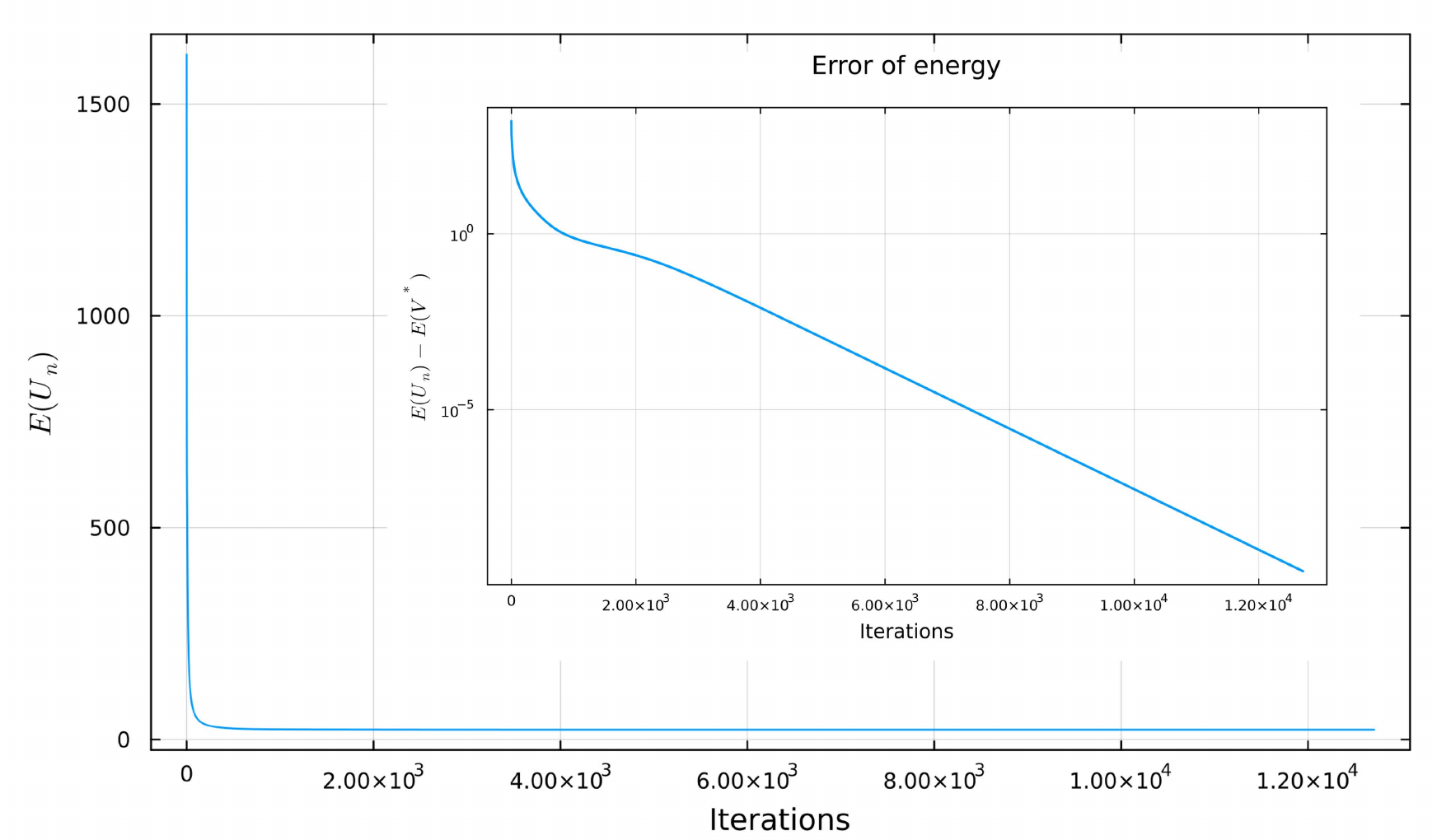}
					\caption{Convergence curves of the energy}
					\label{fig:laplaceenergy}
				\end{subfigure}
				\hfill %
				\begin{subfigure}[b]{0.48\linewidth}
					\includegraphics[width=\linewidth]{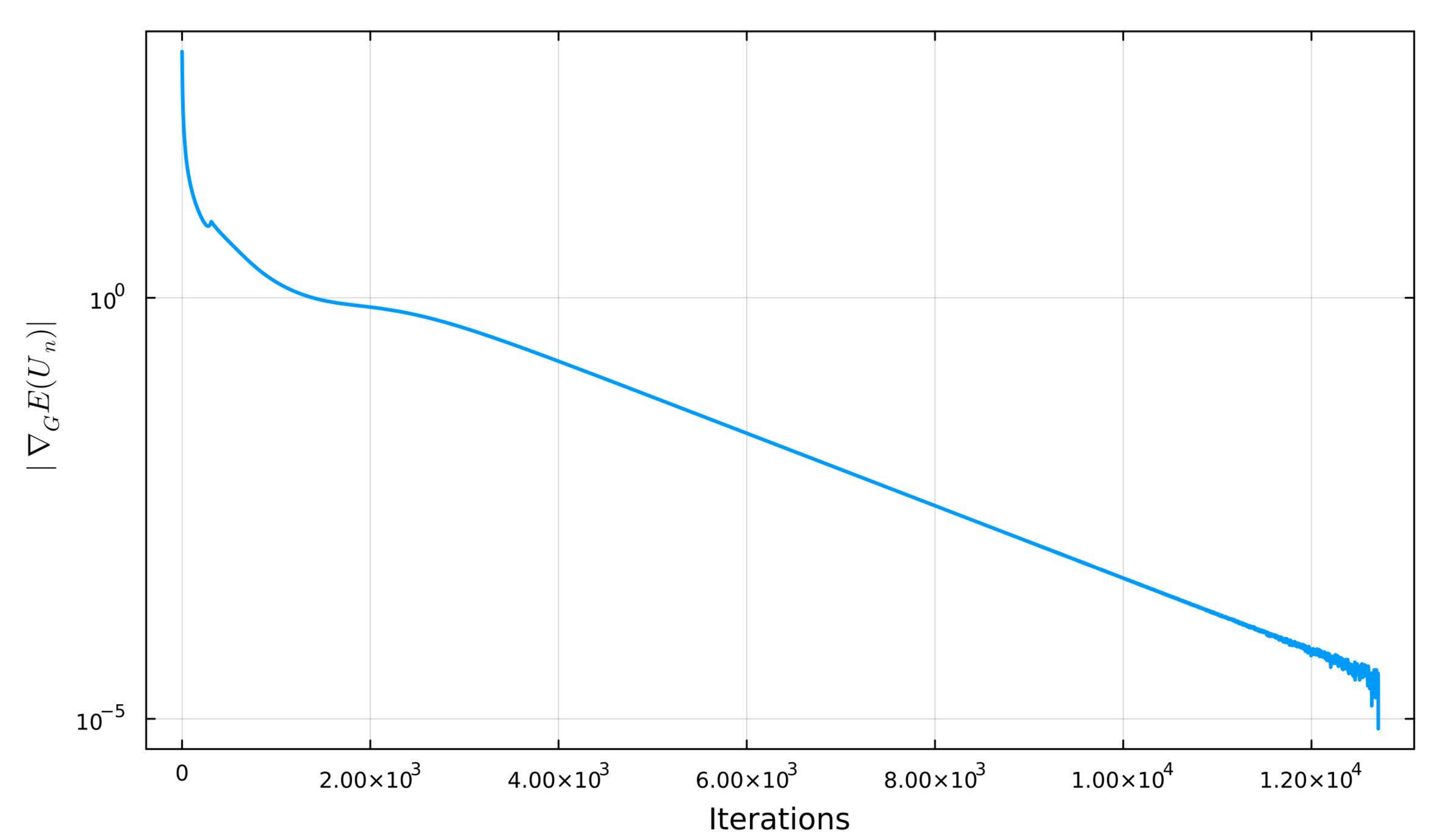}
					\caption{Convergence curve of the gradient}
					\label{fig:laplacegradient}
				\end{subfigure}\\[-3mm] 
				\begin{subfigure}[b]{0.48\linewidth}
					\includegraphics[width=\linewidth]{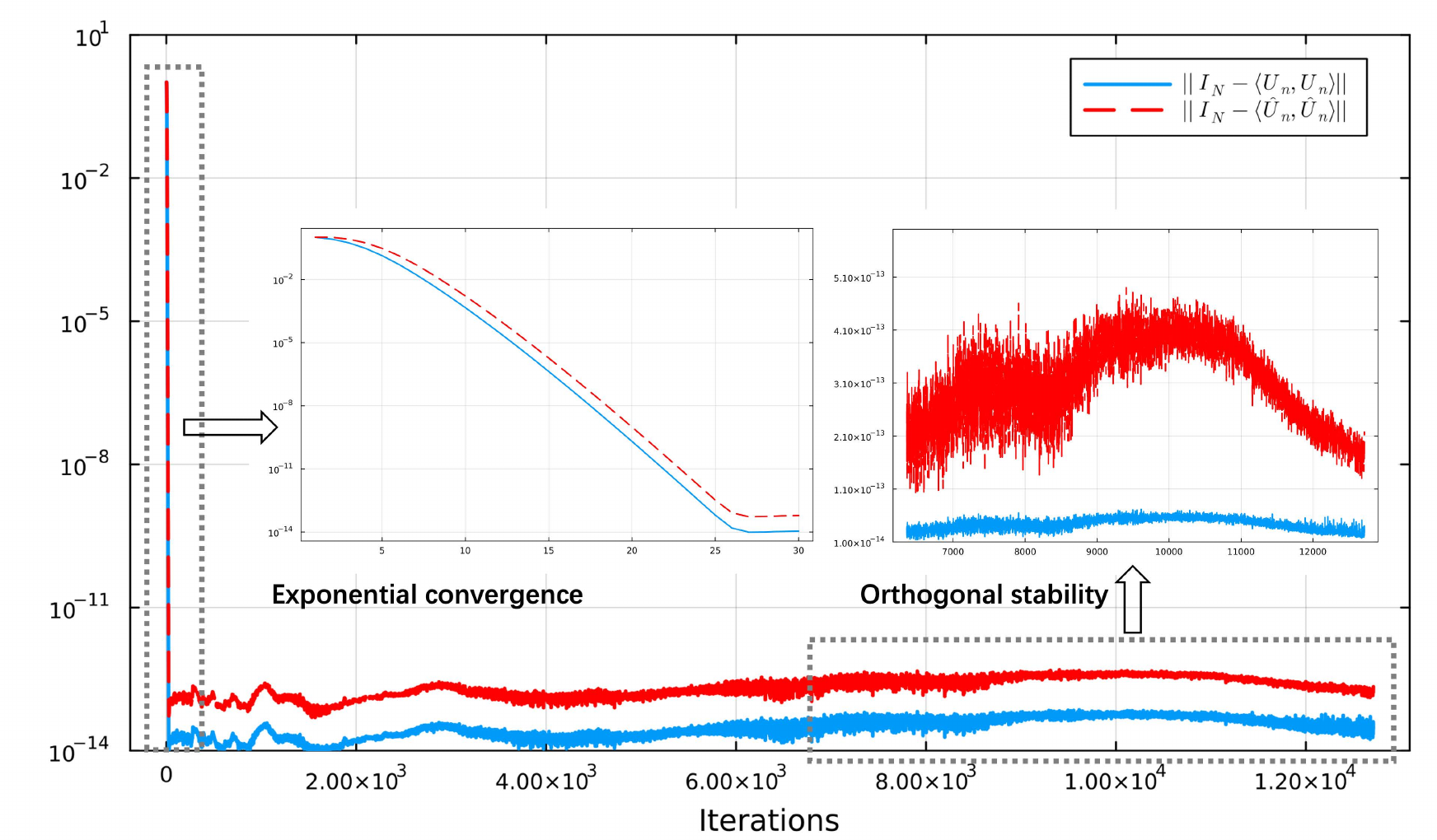}
					\caption{Convergence of the orthogonality error}
					\label{fig:laplaceortho}
				\end{subfigure}
				\hfill
				\begin{subfigure}[b]{0.48\linewidth}
					\includegraphics[width=\linewidth]{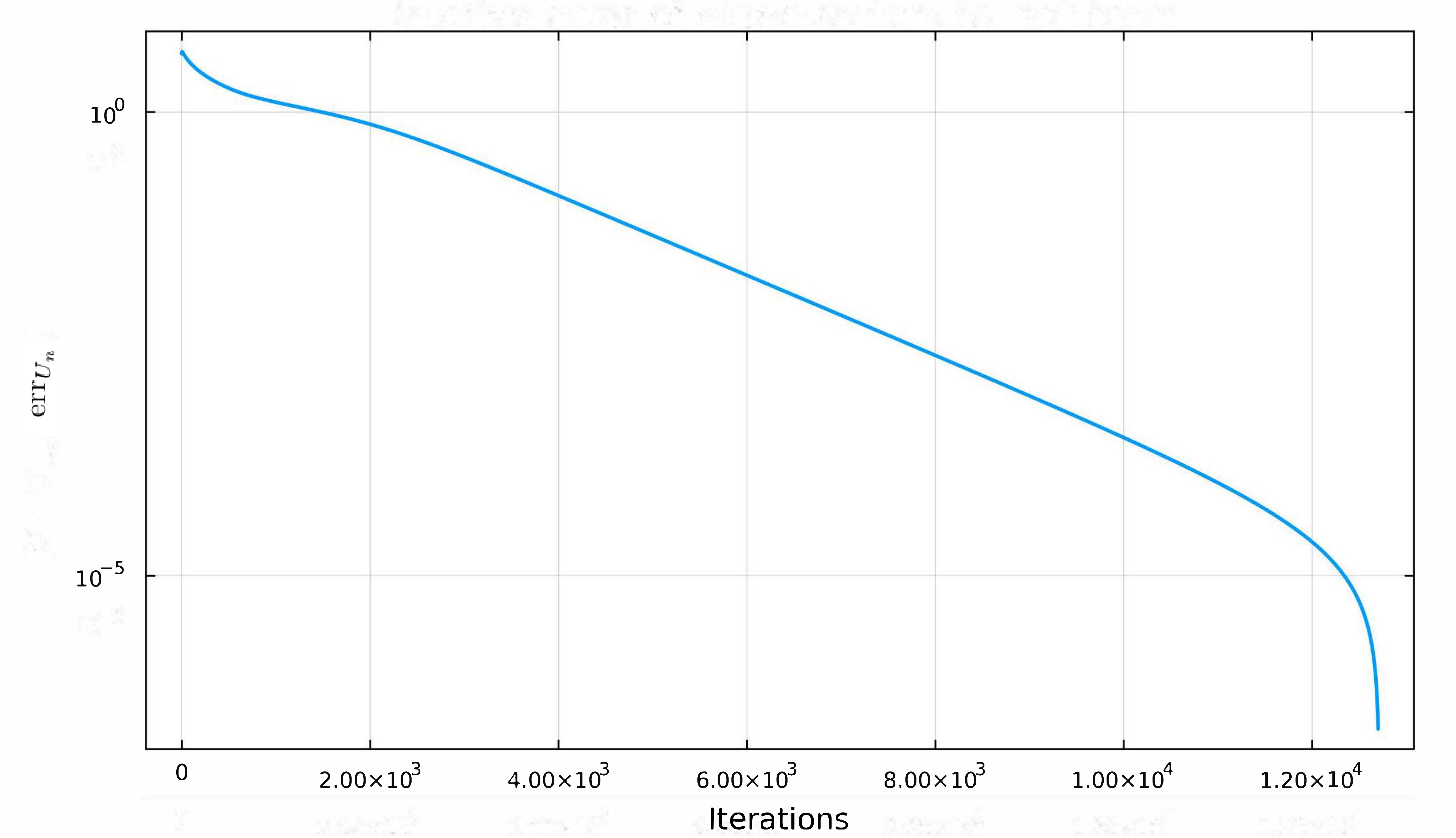}
					\caption{Relative error curve of the iterates}
					\label{fig:laplaceeigenfunerr}
				\end{subfigure}\\[-3mm] 
				\caption{Numerical results of Example~\ref{eq:3D Laplace eigenvalue equation}}
				\label{fig:laplace_total}
			\end{figure}
		\begin{itemize}
			\item Fig.~\ref{fig:laplaceenergy}: The computed energy decreases over the iterations and approaches the reference value used in the plot.
			\item Fig.~\ref{fig:laplacegradient}: The computed gradient norm \(\|\nabla_G E(U_n)\|\) decreases and reaches the prescribed tolerance. At the final iterate, it is also the residual norm of the corresponding eigenpair approximations.
			\item Fig.~\ref{fig:laplaceortho}: The orthogonality error \(\|I_N - \langle U_n,U_n\rangle\|\) decreases during the early iterations and thereafter stabilizes with minor fluctuations. The corrector step (blue line) produces smaller orthogonality errors than the predictor step (red line), and the error reaches the prescribed tolerance.
			\item Fig.~\ref{fig:laplaceeigenfunerr} illustrates the variation of the relative quantity \(\text{err}_{U_n}\) defined in (\ref{eq:relative error of iterates}). This quantity decreases as the iterates approach the final computed iterate.
			\end{itemize}

		\end{example}

	To examine the behavior of the algorithm for different self-adjoint operators, we further consider two eigenvalue problems governed by the harmonic oscillator and hydrogen atom Schr\"odinger operators.
		\begin{example}\label{eq:3D harmonic oscillator equation 3D}
			We consider the harmonic oscillator eigenvalue problem on $\mathbb{R}^3$ \cite{ReedSimonIV}:
			\begin{equation*}
				-\frac{1}{2} \Delta u + \frac{1}{2}|x|^2 u = \lambda u, \qquad \int_{\mathbb{R}^3} u^2 = 1,
			\end{equation*}
			where $(u, \lambda) \in H_0^1(\mathbb{R}^3) \times \mathbb{R}$ and $|x| = \sqrt{x_1^2 + x_2^2 + x_3^2}$. Explicit eigenvalues and eigenfunctions (via Hermite polynomials $\mathcal{H}_n$) are:
			\[
			\lambda_{n_1, n_2, n_3} = n_1 + n_2 + n_3 + \frac{3}{2}, \quad u_{n_1, n_2, n_3}(x) = \mathcal{H}_{n_1}(x_1)\mathcal{H}_{n_2}(x_2)\mathcal{H}_{n_3}(x_3) e^{-\frac{1}{2}|x|^2},
			\]
			for $n_1, n_2, n_3 = 0, 1, \cdots$. Due to exponential decay of eigenfunctions, we restrict to $\Omega = (-5.5, 5.5)^3$ for computation. {We choose \(N=10\) to include the first three distinct eigenvalues with their full multiplicities, so that \(\lambda_N<\lambda_{N+1}\).} A uniform finite element mesh with $N_g = 24389$ degrees of freedom is adopted.
			\begin{figure}[htbp]
				\centering
				\begin{subfigure}[b]{0.48\linewidth}
					\includegraphics[width=\linewidth]{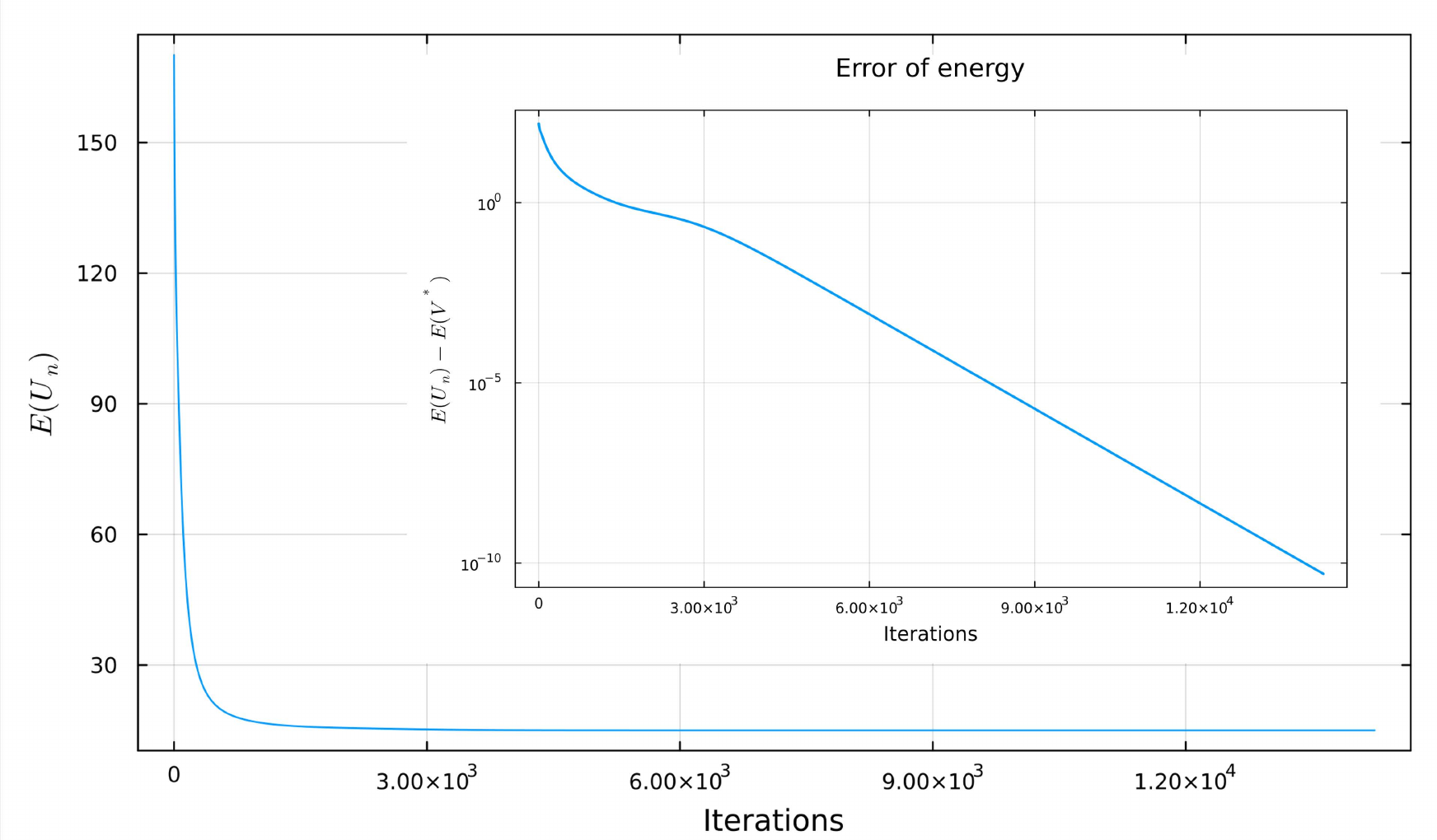}
					\caption{Convergence curves of the energy}
					\label{fig:harmonicenergy}
				\end{subfigure}
				\hfill 
				\begin{subfigure}[b]{0.48\linewidth}
					\includegraphics[width=\linewidth]{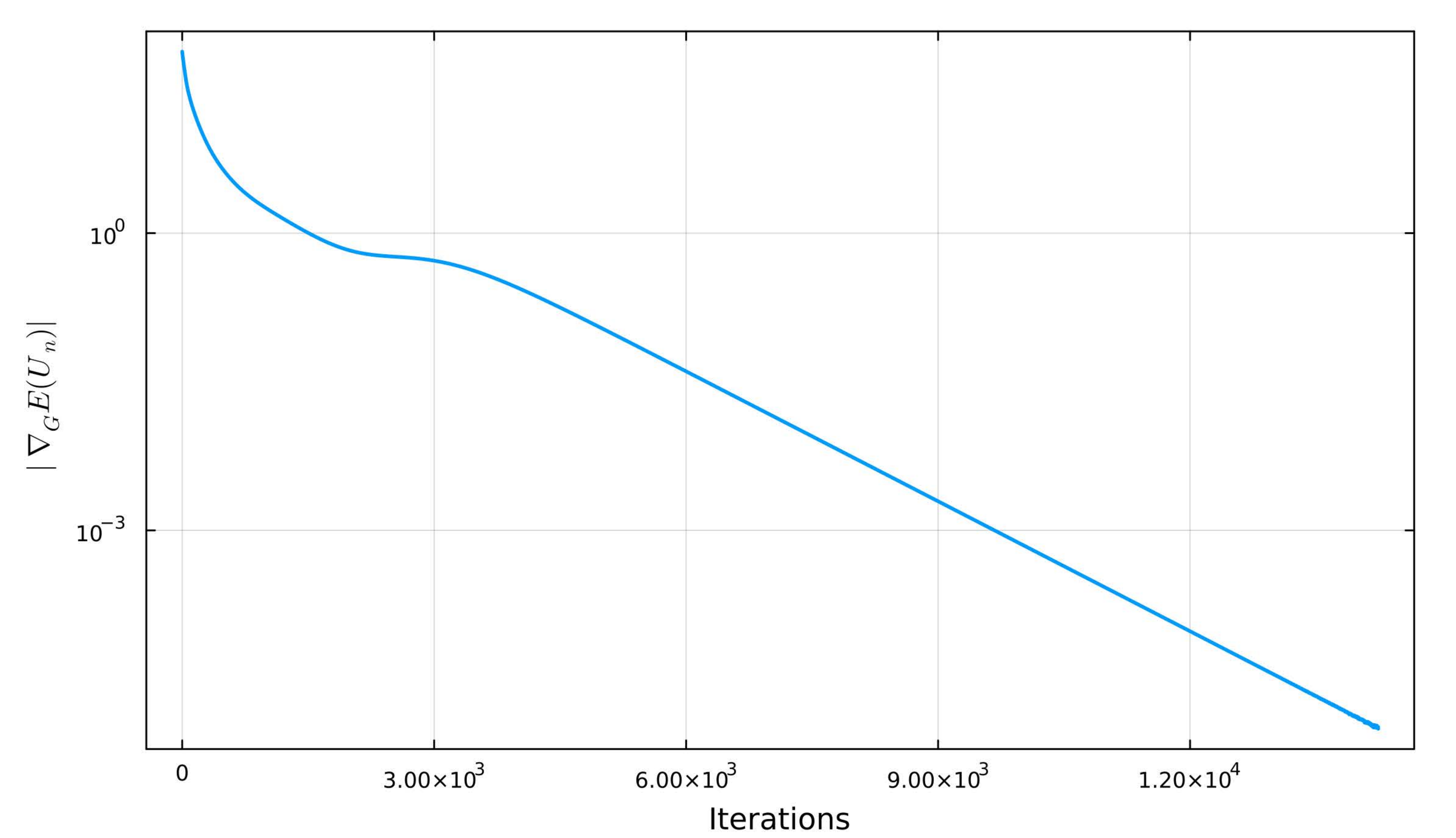}
					\caption{Convergence curve of the gradient}
					\label{fig:harmonicgradient}
				\end{subfigure}\\[-3mm] 
				\begin{subfigure}[b]{0.48\linewidth}
					\includegraphics[width=\linewidth]{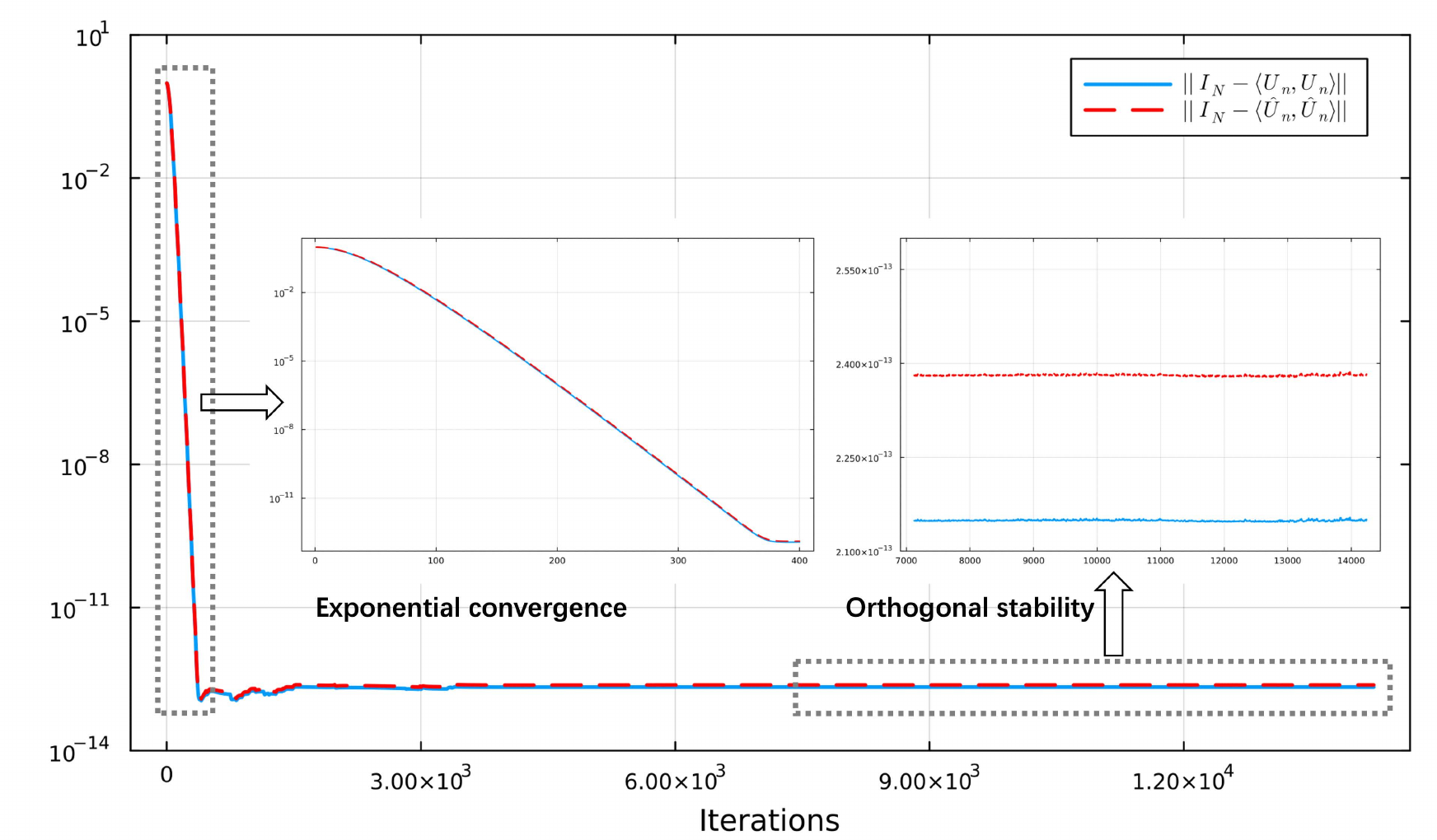}
					\caption{Convergence of the orthogonality error}
					\label{fig:harmonicortho}
				\end{subfigure}
				\hfill
				\begin{subfigure}[b]{0.48\linewidth}
					\includegraphics[width=\linewidth]{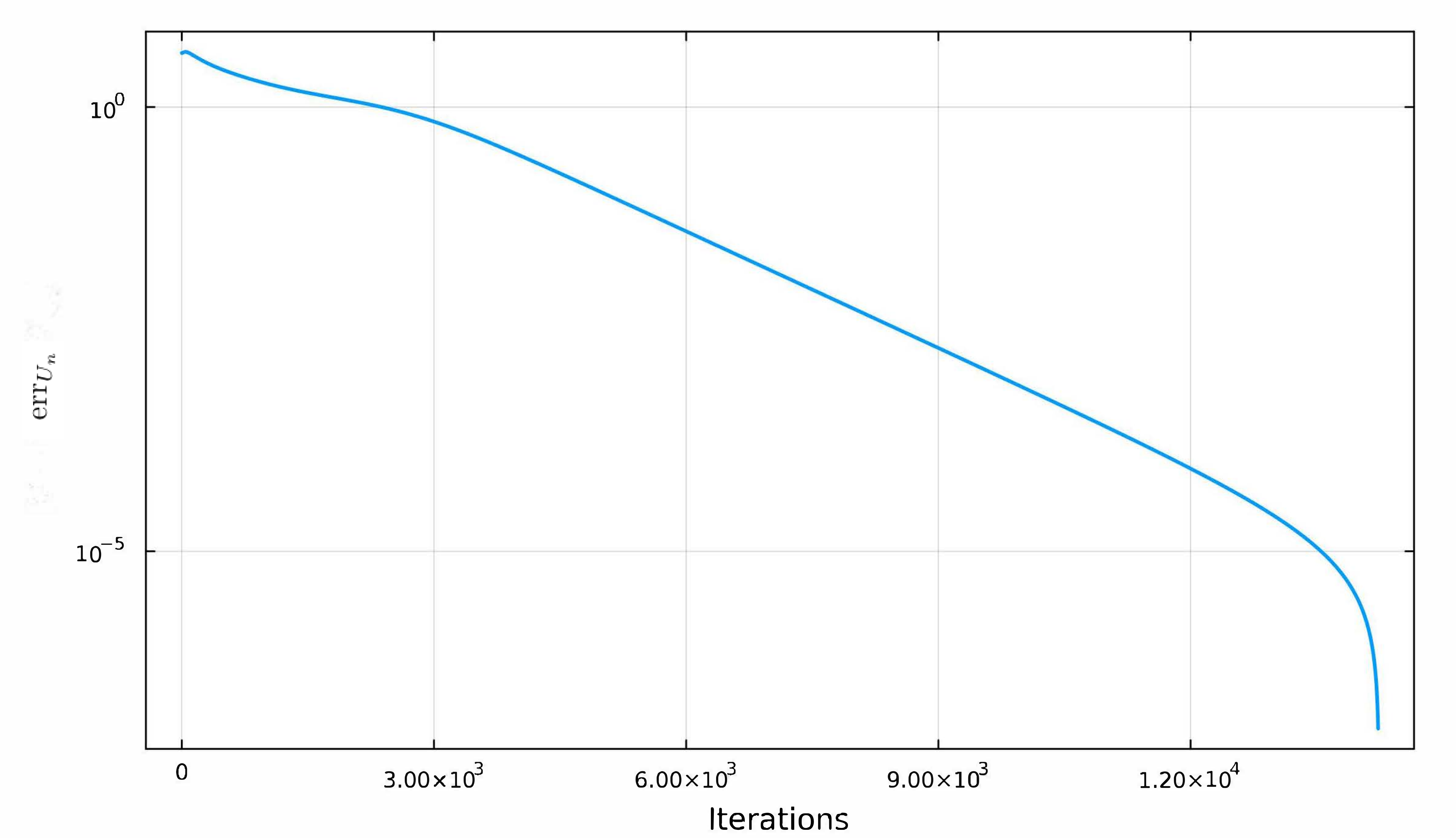}
					\caption{Relative error curve of the iterates}
					\label{fig:harmoniceigenfunerr}
				\end{subfigure}\\[-3mm] 
				\caption{Numerical results of Example~\ref{eq:3D harmonic oscillator equation 3D}}
				\label{fig:harmonic_total}
			\end{figure}
		\end{example}
		
		\begin{example}\label{eq:3D hydrogen}
			We consider the hydrogen atom Schr\"odinger equation on $\mathbb{R}^3$ \cite{greiner2011quantum}:
			\begin{equation*}
				\left(-\frac{1}{2} \Delta - \frac{1}{|x|}\right) u = \lambda u, \qquad \int_{\mathbb{R}^3} |u|^2 = 1,
			\end{equation*}
			where $(u, \lambda) \in H_0^1(\mathbb{R}^3) \times \mathbb{R}$. Eigenvalues are $\lambda_n = -\frac{1}{2n^2}$ ($n=1,2,\cdots$) with multiplicity $n^2$. Using exponential decay of eigenfunctions, we compute on $\Omega = (-20.0, 20.0)^3$. {We choose \(N=5\) to include the first two distinct eigenvalues with their full multiplicities, so that \(\lambda_N<\lambda_{N+1}\).} Spatial discretization uses adaptive finite elements \cite{dai2015convergence}, resulting in $N_g = 13431$ degrees of freedom.
		\begin{figure}[htbp]
				\centering
					\begin{subfigure}[b]{0.48\linewidth}
					\includegraphics[width=\linewidth]{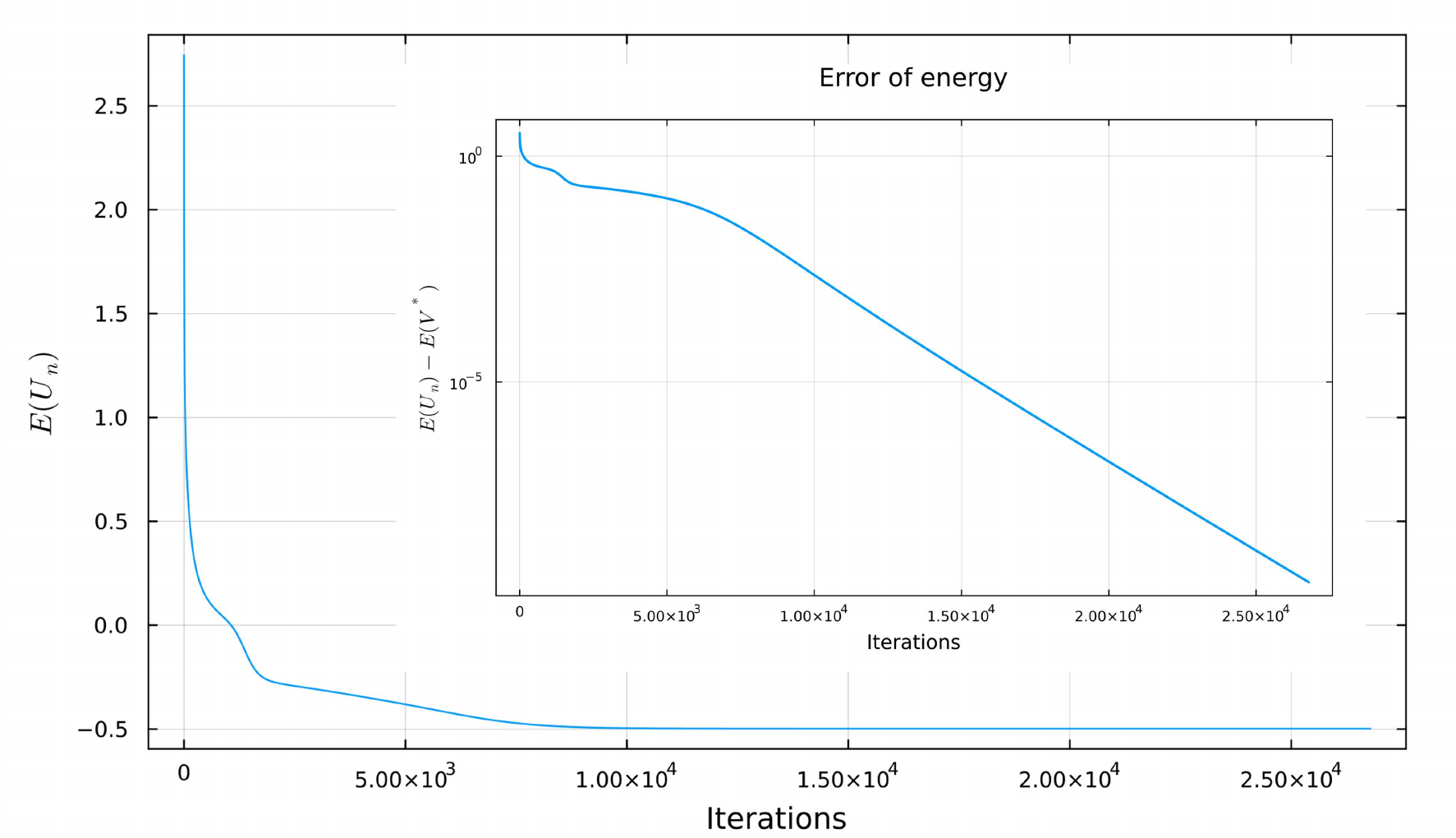}
					\caption{Convergence curves of the energy}
					\label{fig:hydrogenenergy}
				\end{subfigure}
				\hfill 
				\begin{subfigure}[b]{0.48\linewidth}
					\includegraphics[width=\linewidth]{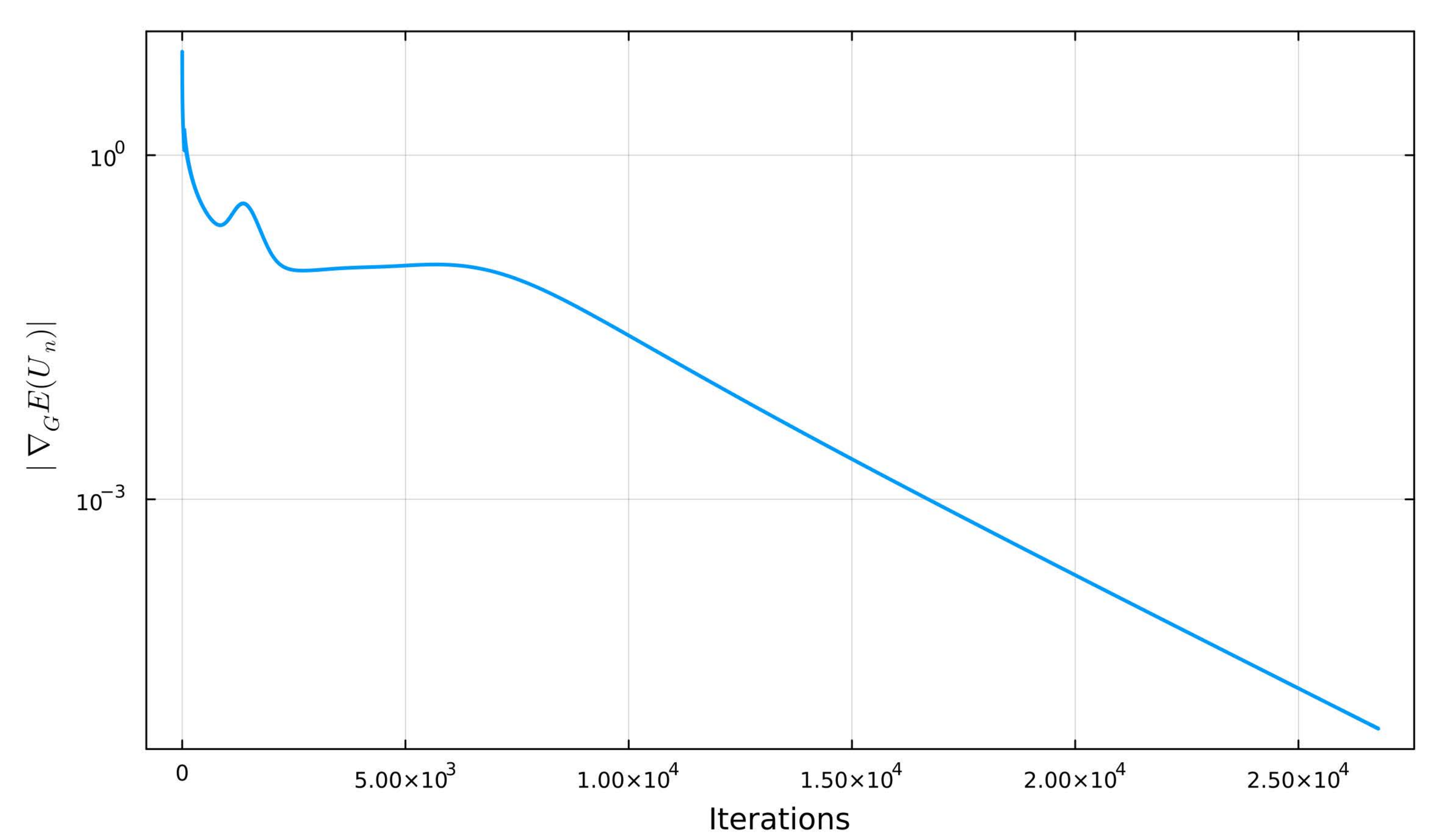}
					\caption{Convergence curve of the gradient}
					\label{fig:hydrogengradient}
				\end{subfigure}\\[-3mm] 
				\begin{subfigure}[b]{0.48\linewidth}
					\includegraphics[width=\linewidth]{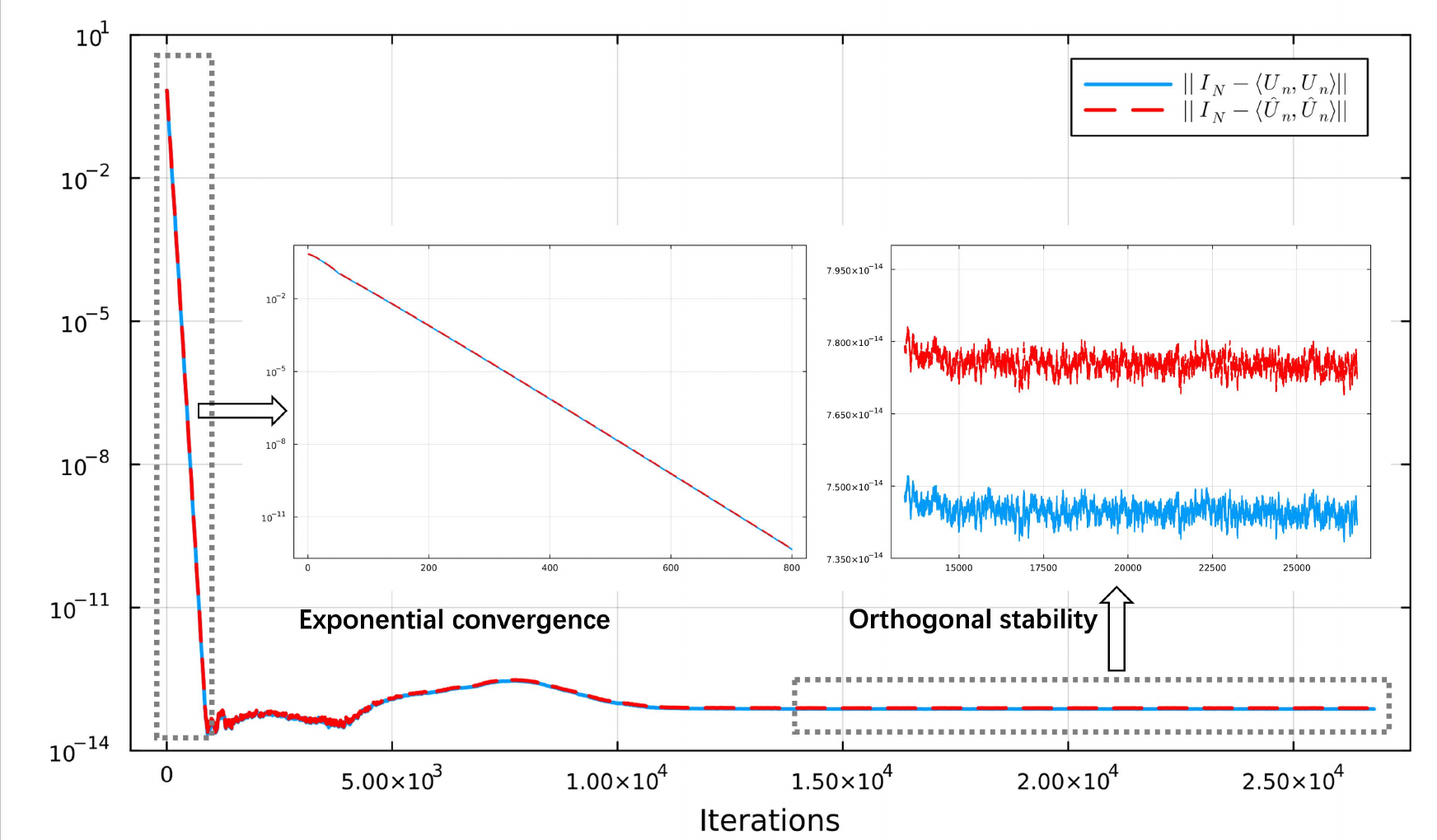}
					\caption{Convergence of the orthogonality error}
					\label{fig:hydrogenortho}
				\end{subfigure}
				\hfill
				\begin{subfigure}[b]{0.48\linewidth}
					\includegraphics[width=\linewidth]{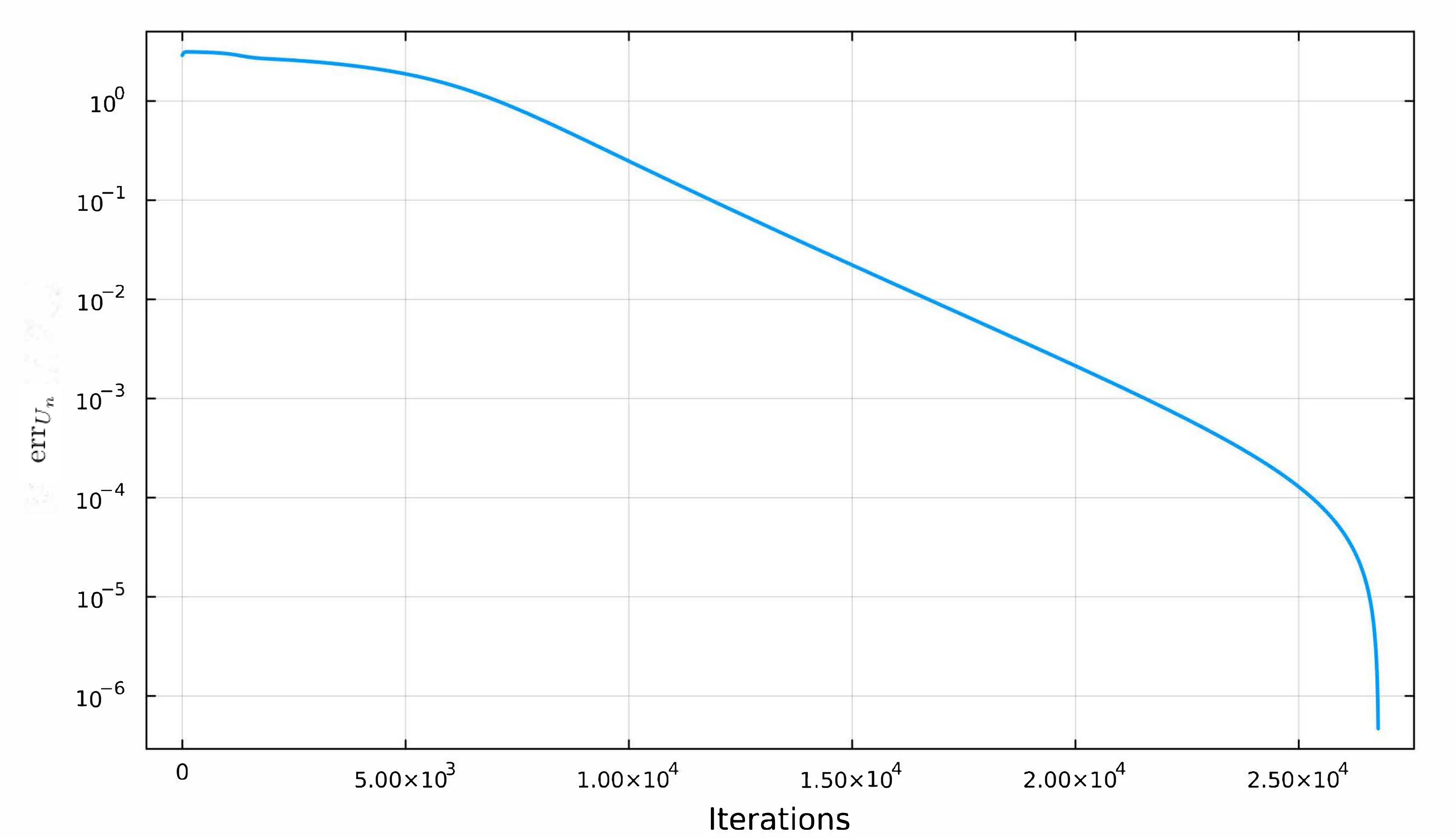}
					\caption{Relative error curve of the iterates}
					\label{fig:hydrogeneigenfunerr}
				\end{subfigure}\\[-3mm] 
				\caption{Numerical results of Example~\ref{eq:3D hydrogen}}
				\label{fig:hydrogen_total}
			\end{figure}
		\end{example}
		
Conclusions for Examples~\ref{eq:3D harmonic oscillator equation 3D} and \ref{eq:3D hydrogen} are qualitatively consistent with that of Example~\ref{eq:3D Laplace eigenvalue equation}, with their numerical results (Figs.~\ref{fig:harmonic_total} and \ref{fig:hydrogen_total}) adopting the same subfigure structure as Fig.~\ref{fig:laplace_total} to eliminate redundant explanations:
\begin{itemize}
	\item Figs.~\ref{fig:harmonicenergy} and \ref{fig:hydrogenenergy}: The computed energies decrease over the iterations, as in the Laplace case.
	\item Figs.~\ref{fig:harmonicgradient} and \ref{fig:hydrogengradient}: The computed gradient norms \(\|\nabla_G E(U_n)\|\) decrease and reach the prescribed tolerance. At the final iterates, they are also the residual norms of the corresponding eigenpair approximations.
	\item Figs.~\ref{fig:harmonicortho} and \ref{fig:hydrogenortho}: The corrector step \(U_n\) (blue) produces smaller orthogonality errors than the predictor step \(\hat{U}_n\) (red). The errors decrease during the early iterations and then stabilize near the prescribed tolerances.
	\item Figs.~\ref{fig:harmoniceigenfunerr} and \ref{fig:hydrogeneigenfunerr}: The relative quantities decrease as the iterates approach their respective final computed iterates.
\end{itemize}

	\section{Concluding remarks} \label{sec:Concluding remarks}
	{We develop and analyze a predictor--corrector discretization of the quasi-Grassmannian gradient flow in \cite{Wang2025} for the simultaneous approximation of a cluster of eigenpairs.} Its iterates approach orthogonality without any orthogonalization operation. For the nonorthogonal initial data considered in the analysis, the iteration preserves full column rank and converges to the target eigenspace from a neighborhood of the solution. {In the numerical experiments, the equations in \eqref{Discretization scheme} are solved approximately, and the computed iterates attain numerical orthogonality to high accuracy, indicating the potential of the proposed scheme for simultaneously computing a cluster of eigenpairs.}

	{The proposed discretization provides a basis for further numerical developments of the quasi-Grassmannian gradient flow. In particular, alternative discretizations and acceleration strategies deserve further investigation, together with extensions to broader operator classes or nonlinear eigenvalue problems and parallel implementations.}

	{\section*{Acknowledgments}
	The authors thank Dr. Liwei Zhang for his careful review of the manuscript, which helped improve the presentation of this paper.}
		
	\appendix
	\section{Proofs of auxiliary results}
	\label{app:auxiliary-proofs}

	\begin{proof}[Proof of Lemma \ref{lemma: existence of g}]
		For the predictor equation, define
		\begin{equation*}
			\Phi(\hat{W},U,\tau)
			=
			\hat{W}-U+\tau\mathcal{A}_{\frac{\hat{W}+U}{2}}
			\frac{\hat{W}+U}{2}.
		\end{equation*}
		For every $U\in[V_*]$,
		\begin{equation*}
			\Phi(U,U,0)=0,
			\qquad
			D_{\hat{W}}\Phi(U,U,0)=\mathcal{I}.
		\end{equation*}
		The implicit function theorem gives a unique solution $\hat{W}=\hat{g}(U,\tau)$ for $U$ near $[V_*]$ and for sufficiently small $\tau$.

		For the corrector equation, define
		\begin{equation*}
			\Psi(W,\hat{W},\tau)
			=
			W-\hat{W}+\tau\nabla E(W)(I_N-\langle W,W\rangle).
		\end{equation*}
		For every $\hat{W}\in[V_*]$,
		\begin{equation*}
			\Psi(\hat{W},\hat{W},0)=0,
			\qquad
			D_W\Psi(\hat{W},\hat{W},0)=\mathcal{I}.
		\end{equation*}
		A second application of the implicit function theorem gives a unique solution $W=g(U,\tau)$ after setting $\hat{W}=\hat{g}(U,\tau)$. {The equations in \eqref{Discretization scheme} are equivariant under right multiplication by an orthogonal matrix.} In particular,
		\begin{equation*}
			\Phi(\hat{W}Q,UQ,\tau)=\Phi(\hat{W},U,\tau)Q,
			\qquad
			\Psi(WQ,\hat{W}Q,\tau)=\Psi(W,\hat{W},\tau)Q.
		\end{equation*}
		The orthogonal invariance of $\hat{g}$ and $g$ follows from uniqueness in this neighborhood. The proof is complete.
	\end{proof}

	\begin{proof}[Proof of Lemma \ref{lemma of energy decreases}]
		Since $ U \in \left(\mathcal{V}_0\right)^N\bigcap \mathcal{M}_{\leqslant}^{N;N_g}$, Theorem \ref{Un inside Stifel manifold} implies
		\begin{equation*} 
				\hat{g}(U, \tau)\in \left(\mathcal{V}_0\right)^N\bigcap \mathcal{M}_{\leqslant}^{N;N_g} \text{ and }	g(U, \tau)\in \left(\mathcal{V}_0\right)^N\bigcap \mathcal{M}_{\leqslant}^{N;N_g}  \quad \forall \tau \in [0, \delta_e]. 
		\end{equation*}
		By the corrector equation,
		\begin{equation*}
			\begin{aligned}
				\hat{g}(U,\tau)=g(U,\tau)+\tau\nabla E(g(U,\tau))\left(I_N-\langle g(U,\tau),g(U,\tau)\rangle\right).
			\end{aligned}
		\end{equation*}
		Since \(\nabla E(g(U,\tau))=Hg(U,\tau)\), we have
		\begin{equation*}
			\begin{aligned}
				&E(\hat{g}(U,\tau))-E(g(U,\tau))
				\\=&\tau\operatorname{tr}\Big(\left(I_N-\langle g(U,\tau),g(U,\tau)\rangle\right)\langle g(U,\tau),H^2g(U,\tau)\rangle\Big)
				\\&+\frac{\tau^2}{2}\operatorname{tr}\Big(\left(I_N-\langle g(U,\tau),g(U,\tau)\rangle\right)\langle g(U,\tau),H^3g(U,\tau)\rangle
				\left(I_N-\langle g(U,\tau),g(U,\tau)\rangle\right)\Big).
			\end{aligned}
		\end{equation*}
		Because \(g(U,\tau)\in \left(\mathcal{V}_0\right)^N\bigcap \mathcal{M}_{\leqslant}^{N;N_g}\), we have \(0\leqslant I_N-\langle g(U,\tau),g(U,\tau)\rangle\leqslant I_N\). Hence
		\begin{equation*}
			\begin{aligned}
				&\left|\operatorname{tr}\Big(\left(I_N-\langle g(U,\tau),g(U,\tau)\rangle\right)\langle g(U,\tau),H^3g(U,\tau)\rangle
				\left(I_N-\langle g(U,\tau),g(U,\tau)\rangle\right)\Big)\right|
				\\&\leqslant |\lambda_1|
				\operatorname{tr}\Big(\left(I_N-\langle g(U,\tau),g(U,\tau)\rangle\right)\langle g(U,\tau),H^2g(U,\tau)\rangle
				\left(I_N-\langle g(U,\tau),g(U,\tau)\rangle\right)\Big)
				\\&\leqslant |\lambda_1|
				\operatorname{tr}\Big(\left(I_N-\langle g(U,\tau),g(U,\tau)\rangle\right)\langle g(U,\tau),H^2g(U,\tau)\rangle\Big).
			\end{aligned}
		\end{equation*}
		Thus, for \(\tau\in[0,\delta_e]\),
		\begin{equation*}
			\begin{aligned}
				E(\hat{g}(U,\tau))-E(g(U,\tau))
				&\geqslant
				\tau\left(1-\frac{\tau}{2}|\lambda_1|\right)
				\\&\quad\cdot
				\operatorname{tr}\Big(\left(I_N-\langle g(U,\tau),g(U,\tau)\rangle\right)\langle g(U,\tau),H^2g(U,\tau)\rangle\Big)
				\geqslant0.
			\end{aligned}
		\end{equation*}
		
		Next, we analyze \(E(\hat{g}(U, \tau)) - E(U)\). Similarly, denote $\tilde{S}(t) = t \hat{g}(U, \tau) +(1-t)U$ for $t\in [0,1]$, and there exists $\tilde{\xi} \in (0,1)$ such that
		\begin{equation*}
			\begin{aligned}
				&	E(\hat{g}(U, \tau))  - E(U) 
				\\=& \operatorname{tr}\left(\left\langle\nabla E(\tilde{S}(\tilde{\xi})),\hat{g}(U,\tau)-U\right\rangle\right)  
				= -\tau \operatorname{tr}\left(\left\langle\nabla E(\tilde{S}(\tilde{\xi})),\mathcal{A}_{\tilde{U}_{\frac{1}{2}}}\tilde{U}_{\frac{1}{2}}\right\rangle\right)
				\\ =& -\tau \underbrace{ \operatorname{tr}\left(\left\langle\nabla E(\tilde{U}_{\frac{1}{2}}),\mathcal{A}_{\tilde{U}_{\frac{1}{2}}}\tilde{U}_{\frac{1}{2}}\right\rangle\right) }_{\stackrel{\Delta}{=}\mathbf{I_3}}
				+\tau\underbrace{ \operatorname{tr}\left(\left\langle\nabla E(\tilde{U}_{\frac{1}{2}})- \nabla E(\tilde{S}(\tilde{\xi})),\mathcal{A}_{\tilde{U}_{\frac{1}{2}}}\tilde{U}_{\frac{1}{2}}\right\rangle\right)}_{\stackrel{\Delta}{=}\mathbf{I_4}},
			\end{aligned}
		\end{equation*}
		where $\tilde{U}_{\frac{1}{2}} = \frac{\hat{g}(U , \tau)+U }{2}$. It follows from $\lambda(\langle\tilde{U}_{\frac{1}{2}},\tilde{U}_{\frac{1}{2}}\rangle)\subset[0,1]$ that
		\begin{equation*} 
				\left\|\mathcal{A}_{\tilde{U}_{\frac{1}{2}}}\tilde{U}_{\frac{1}{2}}\right\|^2 
				\leqslant \left\|\mathcal{A}_{\tilde{U}_{\frac{1}{2}}}\right\|^2=-\operatorname{tr}\left(\left(\mathcal{A}_{\tilde{U}_{\frac{1}{2}}}\right)^2  \right)
				=2\mathbf{I_3}.
		\end{equation*}
		Since \(\mathcal{A}_{\tilde{U}_{\frac{1}{2}}}\) is skew-symmetric,
		\begin{equation*}
			\begin{aligned}
				\mathbf{I_4}	\leqslant& \left\|\nabla E(\tilde{U}_{\frac{1}{2}} )- \nabla E(\tilde{S}(\tilde{\xi}))\right\|  \left\| \mathcal{A}_{\tilde{U}_{\frac{1}{2}}}\tilde{U}_{\frac{1}{2}}\right\|
				\leqslant |\lambda_1| \left\|\tilde{U}_{\frac{1}{2}}  -\tilde{S}(\tilde{\xi})\right\|\left\| \mathcal{A}_{\tilde{U}_{\frac{1}{2}}}\tilde{U}_{\frac{1}{2}}\right\|
				\\=& |\lambda_1| \left|\tilde{\xi} - \frac{1}{2}\right| \left\|\hat{g}(U, \tau) -U\right\|\left\| \mathcal{A}_{\tilde{U}_{\frac{1}{2}}}\tilde{U}_{\frac{1}{2}}\right\|
				\leqslant \frac{\tau}{2}|\lambda_1| \left\|\mathcal{A}_{\tilde{U}_{\frac{1}{2}}}\tilde{U}_{\frac{1}{2}}\right\|^2.
			\end{aligned}
		\end{equation*}
		Combining the above inequalities and using
		\begin{equation*}
			\begin{aligned}
				E(U) - E(g(U, \tau)) &= (E(U) - E(\hat{g}(U, \tau))) + (E(\hat{g}(U, \tau)) - E(g(U, \tau))) \\&\geqslant E(U) - E(\hat{g}(U, \tau)) = \tau(	\mathbf{I_3}-	\mathbf{I_4}),
			\end{aligned}
		\end{equation*}
		we finally obtain
		\begin{equation*} 
				E(U)-E(g(U, \tau))  \geqslant \tau \left(\frac{1}{2} - \frac{\tau}{2}|\lambda_1|\right)   \left\|\mathcal{A}_{\tilde{U}_{\frac{1}{2}}}\tilde{U}_{\frac{1}{2}}\right\|^2  ,  
		\end{equation*}
		and reach the conclusion.
	\end{proof}

	\begin{proof}[Proof of Lemma \ref{extended gradient is Lip}]
		Let $U_i,U_j\in B\left(V_*,\max\{\eta_a,\eta_b\}\right)$. By the definition of $\alpha$,
		\begin{equation*}
			\sigma_{\max}(U_i)\leqslant\alpha
			\quad\text{and}\quad
			\sigma_{\max}(U_j)\leqslant\alpha.
		\end{equation*}
		First,
		\begin{equation*}
			\begin{aligned}
				&\|HU_i\langle U_i,U_i\rangle-HU_j\langle U_j,U_j\rangle\|
				\\&\leqslant \|H\|\|U_i(\langle U_i,U_i\rangle-\langle U_j,U_j\rangle)\|
				+\|H\|\|(U_i-U_j)\langle U_j,U_j\rangle\|
				\\&\leqslant 3\alpha^2\max\{|\lambda_1|,|\lambda_{\max}|\}\|U_i-U_j\|.
			\end{aligned}
		\end{equation*}
		Similarly,
		\begin{equation*}
			\begin{aligned}
				&\|U_i\langle U_i,HU_i\rangle-U_j\langle U_j,HU_j\rangle\|
				\\&\leqslant
				\|U_i(\langle U_i,HU_i\rangle-\langle U_j,HU_j\rangle)\|
				+\|(U_i-U_j)\langle U_j,HU_j\rangle\|
				\\&\leqslant 3\alpha^2\max\{|\lambda_1|,|\lambda_{\max}|\}\|U_i-U_j\|.
			\end{aligned}
		\end{equation*}
		Combining these two estimates gives
		\begin{equation*}
			\begin{aligned}
				\|\mathcal{A}_{U_i}U_i-\mathcal{A}_{U_j}U_j\|
				\leqslant L\|U_i-U_j\|.
			\end{aligned}
		\end{equation*}
		For the second estimate, let $W\in\left(\mathcal{V}^{N_g}\right)^N$. Then
		\begin{equation*}
			\begin{aligned}
				(\mathcal{A}_{U_i}-\mathcal{A}_{U_j})W
				=&\,H(U_i-U_j)\langle U_i,W\rangle
				+HU_j\langle U_i-U_j,W\rangle
				\\&-(U_i-U_j)\langle HU_i,W\rangle
				-U_j\langle H(U_i-U_j),W\rangle.
			\end{aligned}
		\end{equation*}
		Thus
		\begin{equation*}
				\|\mathcal{A}_{U_i}-\mathcal{A}_{U_j}\|
				\leqslant 4\alpha\max\{|\lambda_1|,|\lambda_{\max}|\}\|U_i-U_j\|
			=\hat{L}\|U_i-U_j\|.
		\end{equation*}
		The proof is complete.
	\end{proof}

	\begin{proof}[Proof of Lemma \ref{discretization is bounded}]
		If $U \in B\left([V_{*}], \eta_2\right)$, then there exists a $\tilde{Q} \in \mathcal{O}^N$ such that $\tilde{U} =U\tilde{Q}\in [U]$ and
		$$\|\tilde{U} - V_{*}\|= \operatorname{dist}\left([\tilde{U}], [V_{*}]\right) \leqslant \eta_2.$$
		From the existence in Lemma \ref{lemma: existence of g}, \(\hat{g}(\tilde{U}, \tau) \in B(V_{*}, \eta_b)\) for all \(\tau \in [0, \delta^*]\).
		
		We first estimate \(\operatorname{dist}([\hat{g}(\tilde{U}, \tau)], [V_{*}])\):
		\begin{equation*}
			\begin{aligned}
				&\operatorname{dist}\left([\hat{g}(\tilde{U} ,\tau)] , [V_{*}]\right) 
				= \inf\limits_{Q\in \mathcal{O}^N}\left\|\hat{g}(\tilde{U} ,\tau) - V_{*} Q\right\|
				\\\leqslant& \left\|\tilde{U}  - V_{*} \right\| +\left\|\hat{g}(\tilde{U} ,\tau)- \tilde{U}\right\|
				\leqslant \eta_2 +\tau \left\|\mathcal{A}_{\frac{\hat{g}(\tilde{U} , \tau)+\tilde{U} }{2}}\frac{\hat{g}(\tilde{U} , \tau)+\tilde{U} }{2} \right\|.
			\end{aligned}
		\end{equation*}
		By Lemma \ref{extended gradient is Lip} and \(\frac{\hat{g}(\tilde{U}, \tau)+\tilde{U}}{2} \in B(V_{*}, \max\{\eta_a, \eta_b\})\), we arrive at
		\begin{equation*}
			\begin{aligned}
				\left\|\mathcal{A}_{\frac{\hat{g}(\tilde{U} , \tau)+\tilde{U} }{2}}\frac{\hat{g}(\tilde{U} , \tau)+\tilde{U} }{2}\right\| 
				&= \left\|\mathcal{A}_{\frac{\hat{g}(\tilde{U} , \tau)+\tilde{U} }{2}}\frac{\hat{g}(\tilde{U} , \tau)+\tilde{U} }{2}- \mathcal{A}_{V_*}V_{*}\right\| 
				\\&\leqslant L\left\|\frac{\hat{g}(\tilde{U} , \tau)+\tilde{U} }{2} - V_{*} \right\| \leqslant L\frac{\eta_a+\eta_b}{2} \leqslant L \max\{\eta_a, \eta_b\},
			\end{aligned}
		\end{equation*}
		which implies
		\begin{equation*} 
				\operatorname{dist}\left([\hat{g}(\tilde{U} ,\tau)] , [V_{*}]\right) \leqslant \eta_2 + \tau L \max\{\eta_a, \eta_b\}. 
		\end{equation*}
		From Lemma \ref{lemma: existence of g}, there exists \(Q_\tau\in\mathcal{O}^N\) such that \(\|g(\tilde{U},\tau)Q_\tau-V_*\|\leqslant \eta_b\). Hence
		\begin{equation*}
			\begin{aligned}
				&\left\|\nabla E\left(g(\tilde{U} ,\tau)\right)\left(I_N -\langle g(\tilde{U},\tau),g(\tilde{U},\tau)\rangle \right)\right\|
				\\\leqslant& |\lambda _1| \sigma_{\max}\left(g(\tilde{U} ,\tau)\right)\left\|I_N -\langle g(\tilde{U},\tau),g(\tilde{U},\tau)\rangle \right\|
				\\\leqslant& |\lambda _1| \alpha\left\|\langle V_*,V_*\rangle-Q_\tau^\top\langle g(\tilde{U},\tau),g(\tilde{U},\tau)\rangle Q_\tau\right\|
				\\\leqslant& |\lambda _1| \alpha(\alpha+1)\left\|g(\tilde{U},\tau)Q_\tau-V_*\right\|
				\leqslant |\lambda _1| \alpha(\alpha+1)\eta_b.
			\end{aligned}
		\end{equation*}
		Consequently, we are able to estimate the distance from \([g(\tilde{U}, \tau)]\) to \([V_*]\):
		\begin{equation*}
			\begin{aligned}
				\operatorname{dist}\left([g(\tilde{U} ,\tau)] , [V_{*}]\right) \leqslant& \operatorname{dist}\left([\hat{g}(\tilde{U} ,\tau)] , [V_{*}]\right)  +\left\|g(\tilde{U} ,\tau) -\hat{g}(\tilde{U} ,\tau) \right\|
				\\\leqslant& \eta_2 + \tau L \max\{\eta_a, \eta_b\} 
				\\&+\tau\left\|\nabla E\left(g(\tilde{U} ,\tau)\right)
				\left(I_N -\langle g(\tilde{U},\tau),g(\tilde{U},\tau)\rangle \right)\right\|
				\\ \leqslant &\eta_2 + \tau L \max\{\eta_a, \eta_b\} + \tau |\lambda _1| \alpha(\alpha+1)\eta_b.
			\end{aligned}
		\end{equation*}
		By the orthogonal invariance of $g$, we have
		\begin{equation*}
			\begin{aligned}
			&	\operatorname{dist}\left([g(U ,\tau)] , [V_{*}]\right)  = \operatorname{dist}\left([g(U ,\tau)\tilde{Q}] , [V_{*}]\right) 
				\\= &  \operatorname{dist}\left([g(U\tilde{Q} ,\tau)] , [V_{*}]\right)  
				= \operatorname{dist}\left([g(\tilde{U} ,\tau)] , [V_{*}]\right) \leqslant \eta_1 \quad \forall \tau \in [0, \delta_b] ,
			\end{aligned}
		\end{equation*}
		that is,
		\begin{equation*} 
				g(U,\tau)\in B([V_{*}],\eta_1) \quad \forall \tau \in [0, \delta_b]. 
		\end{equation*}
		As such, the proof is concluded.
	\end{proof}

	\section{Practical implementation details}
	\label{app:practical-implementation}
		The predictor equation can be written as
		\begin{equation}
			\label{splited into two equations}
			\left\{
			\begin{aligned}
				&\frac{\tilde{U}_{n+\frac{1}{2}}-U_n}{\frac{1}{2}\tau_n}
				=-\mathcal{A}_{\tilde{U}_{n+\frac{1}{2}}}\tilde{U}_{n+\frac{1}{2}},
				\\
				&\frac{\hat{U}_{n+1}-\tilde{U}_{n+\frac{1}{2}}}{\frac{1}{2}\tau_n}
				=-\mathcal{A}_{\tilde{U}_{n+\frac{1}{2}}}\tilde{U}_{n+\frac{1}{2}}.
			\end{aligned}
			\right.
		\end{equation}
		The first equation is equivalent to
		\begin{equation}
			\label{eq:midpoint form}
			\tilde{U}_{n+\frac{1}{2}}
			=
			\left(
			\mathcal{I}+\frac{\tau_n}{2}\mathcal{A}_{\tilde{U}_{n+\frac{1}{2}}}
			\right)^{-1}U_n.
		\end{equation}
		Starting from \(\tilde{U}_{n+\frac{1}{2}}^{(0)}=U_n\), we compute
		\begin{equation*}
			\tilde{U}_{n+\frac{1}{2}}^{(k)}
			=
			\left(
			\mathcal{I}+\frac{\tau_n}{2}\mathcal{A}_{\tilde{U}_{n+\frac{1}{2}}^{(k-1)}}
			\right)^{-1}U_n,
			\qquad k=1,\ldots,p_n,
		\end{equation*}
		and set
		\begin{equation*}
			\hat{U}_{n+1}=2\tilde{U}_{n+\frac{1}{2}}^{(p_n)}-U_n.
		\end{equation*}
		We evaluate this inverse using the Sherman--Morrison--Woodbury formula:
		\begin{equation*}
			\begin{aligned}
				\left(\mathcal{I}+\frac{\tau_n}{2}\mathcal{A}_U\right)^{-1}
				={}&\mathcal{I}-\frac{\tau_n}{2}
				\left(\begin{array}{cc}
					\nabla E(U) & U
				\end{array}\right)
				\\
				&\times\left[
				I_{2N}+\frac{\tau_n}{2}
				\left(\begin{array}{cc}
					\langle U,\nabla E(U)\rangle & \langle U,U\rangle
					\\
					-\langle\nabla E(U),\nabla E(U)\rangle & -\langle\nabla E(U),U\rangle
				\end{array}\right)
				\right]^{-1}
				\\
				&\times
				\binom{\langle U,\cdot\rangle}{-\langle\nabla E(U),\cdot\rangle}.
			\end{aligned}
		\end{equation*}
		Thus, only a \(2N\times2N\) matrix is inverted in this computation.

		The corrector equation is
		\begin{equation*}
			U_{n+1}
			=
			\hat{U}_{n+1}
			-\tau_n\nabla E(U_{n+1})
			\left(I_N-\langle U_{n+1},U_{n+1}\rangle\right).
		\end{equation*}
		It may be solved by a nonlinear equation solver to a prescribed residual tolerance. In the numerical experiments, we retain the original implementation and perform one Picard update initialized at \(\hat{U}_{n+1}\):
		\begin{equation*}
			U_{n+1}
			=
			\hat{U}_{n+1}
			-\tau_n\nabla E(\hat{U}_{n+1})
			\left(I_N-\langle\hat{U}_{n+1},\hat{U}_{n+1}\rangle\right).
		\end{equation*}
		Since
		\begin{equation*}
			\begin{aligned}
				\hat{U}_{n+1}
				={}&
				\left(
				\mathcal{I}-\frac{\tau_n}{2}
				\mathcal{A}_{\tilde{U}_{n+\frac{1}{2}}^{(p_n-1)}}
				\right)
				\left(
				\mathcal{I}+\frac{\tau_n}{2}
				\mathcal{A}_{\tilde{U}_{n+\frac{1}{2}}^{(p_n-1)}}
				\right)^{-1}U_n,
			\end{aligned}
		\end{equation*}
		the skew-symmetry of \(\mathcal{A}_{\tilde{U}_{n+\frac{1}{2}}^{(p_n-1)}}\) gives
		\begin{equation*}
			\langle\hat{U}_{n+1},\hat{U}_{n+1}\rangle
			=
			\langle U_n,U_n\rangle.
		\end{equation*}
		Therefore, the implemented corrector becomes
		\begin{equation*}
			U_{n+1}
			=
			\hat{U}_{n+1}
			-\tau_n\nabla E(\hat{U}_{n+1})
			\left(I_N-\langle U_n,U_n\rangle\right).
		\end{equation*}
		This one-step update is the corrector used in Algorithm \ref{alg:Practical Iteration}.

		\bibliographystyle{siamplain}
		\bibliography{references}

@article{Wang2025,
	title={A quasi-{G}rassmannian gradient flow model for eigenvalue problems},
	author={Wang, Shengyue and Zhou, Aihui},
	journal={SIAM J. Math. Anal., to appear},
	note={arXiv:2506.20195}
}

@article{chu2025orthogonality,
	title={An orthogonality-preserving approach for eigenvalue problems},
	author={Chu, Tianyang and Dai, Xiaoying and Wang, Shengyue and Zhou, Aihui},
	journal={arXiv preprint arXiv:2511.06788}
}

@article{zhou2006chebyshev,
	author  = {Zhou, Yunkai and Saad, Yousef and Tiago, Murilo L. and
	           Chelikowsky, James R.},
	title   = {Self-consistent-field calculations using
	           {C}hebyshev-filtered subspace iteration},
	journal = {J. Comput. Phys.},
	volume  = {219},
	number  = {1},
	pages   = {172--184},
	year    = {2006}
}

@article{winkelmann2019chase,
	author  = {Winkelmann, Jan and Springer, Paul and Di Napoli, Edoardo},
	title   = {{ChASE}: A {C}hebyshev accelerated subspace iteration
	           eigensolver for sequences of {H}ermitian eigenvalue problems},
	journal = {ACM Trans. Math. Softw.},
	volume  = {45},
	number  = {2},
	pages   = {21:1--21:34},
	year    = {2019}
}

@article{bathe1973solution,
	title     = {Solution methods for eigenvalue problems in structural mechanics},
	author    = {Bathe, Klaus-J{\"u}rgen and Wilson, Edward L},
	journal   = {Int. J. Numer. Meth. Eng.},
	volume    = {6},
	number    = {2},
	pages     = {213--226},
	year      = {1973},
	publisher = {Wiley Online Library}
}

@article{le2005computational,
	title     = {Computational chemistry from the perspective of numerical analysis},
	author    = {Le Bris, Claude},
	journal   = {Acta Numer.},
	volume    = {14},
	pages     = {363--444},
	year      = {2005},
	publisher = {Cambridge University Press}
}

@book{leszczynski2012handbook,
	title     = {Handbook of Computational Chemistry},
	author    = {Leszczynski, Jerzy},
	volume    = {3},
	year      = {2012},
	publisher = {Springer},
	address   = {Dordrecht}
}

@article{dai2015convergence,
	title     = {Convergence and quasi-optimal complexity of adaptive finite element computations for multiple eigenvalues},
	author    = {Dai, Xiaoying and He, Lianhua and Zhou, Aihui},
	journal   = {IMA J. Numer. Anal.},
	volume    = {35},
	number    = {4},
	pages     = {1934--1977},
	year      = {2015},
	publisher = {Oxford University Press}
}

@article{leon2013gram,
	title     = {Gram-{S}chmidt orthogonalization: 100 years and more},
	author    = {Leon, Steven J and Bj{\"o}rck, {\AA}ke and Gander, Walter},
	journal   = {Numer. Linear Algebra Appl.},
	volume    = {20},
	number    = {3},
	pages     = {492--532},
	year      = {2013},
	publisher = {Wiley Online Library}
}

@article{dai2019adaptive,
	title       = {Adaptive step size strategy for orthogonality constrained line search methods},
	author      = {Dai, Xiaoying and Zhang, Liwei and Zhou, Aihui},
	journal     = {arXiv preprint arXiv:1906.02883}
}

@article{dai2020,
	author  = {Dai, Xiaoying and Wang, Qiao and Zhou, Aihui},
	title   = {Gradient flow based {K}ohn-{S}ham density functional theory model},
	journal = {Multiscale Model. Simul.},
	volume  = {18},
	number  = {4},
	pages   = {1621-1663},
	year    = {2020}
}

@article{dai2021convergent,
	title     = {Convergent and orthogonality preserving schemes for approximating the {K}ohn-{S}ham orbitals},
	author    = {Dai, Xiaoying and Zhang, Liwei and Zhou, Aihui},
	journal   = {Numer. Math. Theor. Meth. Appl.},
	volume    = {16},
	number    = {1},
	pages     = {112204},
	year      = {2023}
}

@article{davidson1975iterative,
	title     = {The Iterative Calculation of a Few of the Lowest Eigenvalues and Corresponding Eigenvectors of Large Real-Symmetric Matrices},
	author    = {E. R. Davidson},
	journal   = {J. Comput. Phys.},
	volume    = {17},
	number    = {1},
	pages     = {87--94},
	year      = {1975},
	publisher = {Elsevier}
}

@article{demmel2012communication,
	title     = {Communication-Optimal Parallel and Sequential {QR} and {LU} Factorizations},
	author    = {Jack Demmel and Laura Grigori and Mark Hoemmen and Julien Langou},
	journal   = {SIAM J. Sci. Comput.},
	volume    = {34},
	number    = {1},
	pages     = {A206--A239},
	year      = {2012},
	publisher = {SIAM}
}

@book{dorfler2011photonic,
	title     = {Photonic Crystals: Mathematical Analysis and Numerical Approximation},
	author    = {D{\"o}rfler, Willy and Lechleiter, Armin and Plum, Michael and Schneider, Guido and Wieners, Christian},
	volume    = {42},
	year      = {2011},
	publisher = {Springer Science \& Business Media},
	address   = {Berlin}
}

@book{evans2022partial,
	author    = {Evans, Lawrence C.},
	title     = {Partial Differential Equations},
	volume    = {19},
	publisher = {American Mathematical Society},
	year      = {2022},
	address   = {Providence}
}

@book{golub2013matrix,
	title     = {Matrix Computations},
	author    = {Gene H. Golub and Charles F. Van Loan},
	year      = {2013},
	edition   = {4th},
	publisher = {Johns Hopkins University Press},
	address   = {Baltimore},
	isbn      = {978-1-4214-0794-4}
}

@book{greiner2011quantum,
	title     = {Quantum Mechanics: An Introduction},
	author    = {Greiner, Walter},
	year      = {2011},
	publisher = {Springer Science \& Business Media},
	address   = {Berlin}
}

@inproceedings{cheng1975eigenfunctions,
	title     = {Eigenfunctions and Eigenvalues of {L}aplacian},
	author    = {Cheng, Shiu Yuen},
	booktitle = {Proc. Symp. Pure Math.},
	volume    = {27},
	number    = {part 2},
	pages     = {185--193},
	year      = {1975}
}

@book{martin2020electronic,
	author    = {Martin, Richard M.},
	title     = {Electronic Structure: Basic Theory and Practical Methods},
	year      = {2020},
	publisher = {Cambridge University Press},
	address   = {Cambridge}
}

@book{ReedSimonIV,
	author    = {Reed, Michael and Simon, Barry},
	title     = {Methods of Modern Mathematical Physics: Analysis of Operators},
	volume    = {4},
	publisher = {Academic Press},
	address   = {New York},
	year      = {1978},
	isbn      = {0-12-585004-2}
}

@book{saad1992numerical,
	title     = {Numerical Methods for Large Eigenvalue Problems},
	author    = {Saad, Youcef},
	year      = {1992},
	publisher = {SIAM},
	address   = {Philadelphia},
	series    = {Classics in Applied Mathematics},
	isbn      = {978-0-89871-534-7}
}

@book{zeidler2013nonlinear,
	title     = {Nonlinear Functional Analysis and Its Applications: II/B: Nonlinear Monotone Operators},
	author    = {Zeidler, Eberhard},
	year      = {2013},
	publisher = {Springer Science \& Business Media},
	address   = {Berlin, Heidelberg},
	series    = {Nonlinear Functional Analysis and Its Applications},
	volume    = {II/B}
}

@article{schneider2009direct,
	title={Direct minimization for calculating invariant subspaces in density functional computations of the electronic structure},
	author={Schneider, Reinhold and Rohwedder, Thorsten and Neelov, Alexey and Blauert, Johannes},
	journal={J. Comput. Math.},
	volume  = {27},
	pages={360--387},
	year={2009},
	publisher={JSTOR}
}
	\end{document}